\documentclass[a4paper,dvipsnames,fleqn]{article}

\usepackage{amsmath}
\usepackage{amssymb}
\usepackage{amsthm}
\usepackage{arydshln}
\usepackage{bm}
\usepackage[colorlinks=true,allcolors=purple]{hyperref}
\usepackage{cleveref}
	\crefformat{equation}{#2(#1)#3}
\usepackage{dsfont}
\usepackage{enumitem}
\usepackage[OT2,OT1]{fontenc}
\usepackage[cal=cm,scr=boondoxo]{mathalfa}
\usepackage{pifont}
\usepackage{stmaryrd}
\usepackage{textcomp}
\usepackage{tikz}
	\usetikzlibrary{arrows}
	\usetikzlibrary{patterns}
\usepackage{tikz-cd}
\usepackage{xfrac}

\numberwithin{equation}{section}			
\swapnumbers						

\newcommand\cyr{
\renewcommand\rmdefault{wncyr}%
\renewcommand\sfdefault{wncyss}%
\renewcommand\encodingdefault{OT2}%
\normalfont
\selectfont}
\DeclareTextFontCommand{\textcyr}{\cyr}

\newcounter{Enum}				
\newenvironment{Enumerate}{\begin{enumerate}[label={\rm({\roman*})}]}{\end{enumerate}}
\newcommand{\Enumref}[1]{{\setcounter{Enum}{#1}{\rm(\roman{Enum})}}}

\newcommand{\descriptionlabelsave}{}		
\newenvironment{Itemize}{%
	\renewcommand{\descriptionlabelsave}{\descriptionlabel}\renewcommand{\descriptionlabel}{$\triangleright$}%
	\begin{description}[leftmargin=15pt,itemindent=-5.2pt]}{%
	\end{description}\renewcommand{\descriptionlabel}{\descriptionlabelsave}}

\newcounter{StepsCount}				

\newcounter{StepsRefCount}

\newenvironment{Ilist}{
	\begin{list}{$\triangleright$}{\leftmargin=0pt \labelwidth=11pt \itemindent=\labelwidth%
	\itemsep=5pt\listparindent=\parindent}}{\end{list}}

\theoremstyle{plain}
	\newtheorem{lemma}{Lemma}[section]
	\newtheorem{proposition}[lemma]{Proposition}
	\newtheorem{theorem}[lemma]{Theorem}
	\newtheorem{corollary}[lemma]{Corollary}
	\newcommand{\GenericTheoremName}{}\newtheorem{generictheorem}[lemma]{\GenericTheoremName}
\theoremstyle{definition}
	\newtheorem{definition}[lemma]{Definition}
	\newcommand{\GenericDefinitionName}{}\newtheorem{genericdefinition}[lemma]{\GenericDefinitionName}
\theoremstyle{remark}
	\newtheorem{remark}[lemma]{Remark}
	\newtheorem{example}[lemma]{Example}
	\newcommand{\GenericRemarkName}{}\newtheorem{genericremark}[lemma]{\GenericRemarkName}
\newenvironment{Proposition}{\begin{proposition}}{\par\noindent\rule{5em}{1pt}\end{proposition}}
\newenvironment{Theorem}{\begin{theorem}}{\par\noindent\rule{5em}{1pt}\end{theorem}}
\newenvironment{Definition}{\begin{definition}}{\par\noindent\rule{5em}{1pt}\end{definition}}
\newenvironment{Remark}{\begin{remark}}{\par\noindent\rule{5em}{0.5pt}\end{remark}}
\newenvironment{Example}{\begin{example}}{\par\noindent\rule{5em}{0.5pt}\end{example}}

\newcommand{\ms}[1]{{\mathscr{#1}}}			
\newcommand{\bb}[1]{{\mathbb{#1}}}			
\newcommand{\ov}{\overline}				

\DeclareMathOperator{\IM}{Im}				
\renewcommand{\Im}{\IM}

\newcommand{\Dis}[1]{${\displaystyle{#1}}$}		
\newcommand{\Side}[1]{\hfill{#1}\kern10pt}		
\newcommand{\FD}[5]{
	\DF\left\{\begin{array}{rcl}{#1}&\to &{#2}\\[#3pt] {#4}&\mapsto &{#5}\end{array}\right.}

\newcommand{\smfrac}[2]{{\textstyle\frac{{#1}}{{#2}}}}	
\newcommand{\smmatrix}[4]{\Bigl(			
	\begin{smallmatrix}
	\hspace*{-0.2ex} #1 \hspace*{0.2ex} & \hspace*{0.2ex} #2 \hspace*{-0.2ex}
	\\[0.5ex]
	\hspace*{-0.2ex} #3 \hspace*{0.2ex} & \hspace*{0.2ex} #4 \hspace*{-0.2ex}
	\end{smallmatrix}
	\Bigr)}

\newcommand{\Dummy}{\text{\textvisiblespace\kern1pt}}	
\newcommand{\Smallo}{{\rm o}}				
\newcommand{\BigO}{{\rm O}}				

\DeclareMathOperator{\Span}{span}			
\DeclareMathOperator{\Ran}{ran}				

\newcommand{\DS}{\mid\mkern3mu}				
\newcommand{\DQ}{\mkern6mu}				
\newcommand{\DP}{{.\kern5pt}}				
\newcommand{\DF}{\colon}				
\newcommand{\DE}{\mathrel{\mathop:}=}			
\newcommand{\ED}{=\mathrel{\mathop:}}			
\newcommand{\DI}{\mathrel{\mathop:}\Leftrightarrow}	
\newcommand{\DD}{\mkern4mu\mathrm{d}}			

\DeclareMathOperator{\Tr}{tr}				

\DeclareMathOperator{\GL}{GL}				

\DeclareMathOperator{\Ind}{Ind}				

\begin{document}

\begin{flushleft}
	{\Large\bf A growth estimate for the monodromy matrix\\[2mm] of a canonical system}
	\\[5mm]
	\textsc{
	Raphael Pruckner
	\,\ $\ast$\,\ 
	Harald Woracek
		\hspace*{-14pt}
		\renewcommand{\thefootnote}{\fnsymbol{footnote}}
		\setcounter{footnote}{2}
		\footnote{This work was supported by the project P~30715-N35 of the 
			Austrian Science Fund (FWF). The second author was supported by the joint project I~4600 of the Austrian
			Science Fund (FWF) and the Russian foundation of basic research (RFBR).}
		\renewcommand{\thefootnote}{\arabic{footnote}}
		\setcounter{footnote}{0}
	}
	\\[6mm]
	{\small
	\textbf{Abstract:}
		We investigate the spectrum of $2$-dimensional canonical systems in the limit circle case. 
		It is discrete and, by the Krein-de~Branges formula, cannot be more dense than the integers. 
		But in many cases it will be more sparse. 
		The spectrum of a particular selfadjoint realisation coincides with the zeroes of one entry of the monodromy 
		matrix of the system. Classical function theory thus establishes an immediate connection between the growth of
		the monodromy matrix and the distribution of the spectrum.

		We prove a generic and flexibel upper estimate for the monodromy matrix, use it to prove a bound for the
		case of a continuous Hamiltonian, and construct examples which show that this bound is sharp. 
		The first two results run along the lines of earlier work of R.Romanov, but significantly improve upon
		these results. This is seen even on the rough scale of exponential order. 
	\\[3mm]
	}
\end{flushleft}

\section{Introduction}

We investigate the spectral theory of $2$-dimensional \emph{canonical systems}
\begin{equation}
\label{Q8}
	y'(t)=zJH(t)y(t),\quad t\in I
	,
\end{equation}
where 
\begin{Itemize}
\item $I=[\alpha,\beta]$ is a finite interval with nonempty interior,
\item $H\DF I\to\bb R^{2\times 2}$ is a (Lebesgue-) measurable function which is integrable and does not vanish on any set of
	positive measure,
\item $H(t)\geq 0$ for almost all $t\in I$,
\item $J$ is the symplectic matrix $J\DE\smmatrix{0}{-1}{1}{0}$,
\item $z$ is a complex parameter (the eigenvalue parameter).
\end{Itemize}
The function $H$ is called the \emph{Hamiltonian} of the system \cref{Q8}. 
Systems of this form are intensively investigated since they can be seen as a unifying framework which includes, e.g., 
Schr\"odinger operators, Jacobi matrices, Dirac operators, and many others. 
Some recent standard literature is \cite{remling:2018,romanov:1408.6022v1,behrndt.hassi.snoo:2020}.

With the system \cref{Q8} one can associate an operator model. It consists of a Hilbert space $L^2(H)$, the maximal and minimal
operators $T_{\max}(H)$ and $T_{\min}(H)$, and a boundary value map $\Gamma(H)\DF T_{\max}(H)\to\bb C^2\times\bb C^2$ (here
we understand $T_{\max}(H)$ as its graph). Selfadjoint realisations of \cref{Q8} have compact resolvents, and are obtained by
specifying boundary conditions on the right and left endpoints $\alpha$ and $\beta$. Each two of them are finite rank
perturbations of each other, and the rank of the perturbation is at most $2$. 
Pick one and denote its eigenvalues as (this sequence need not be two-sided infinite)
\[
	\cdots\leq\lambda_{-2}\leq\lambda_{-1}<0\leq\lambda_1\leq\lambda_2\leq\lambda_3\cdots
\]
By the Krein-de~Branges formula we have (understanding the limit for a finite sequence as $0$)
\begin{equation}
\label{Q9}
	\lim_{n\to\infty}\frac n{\lambda_{|n|}}=\frac 1\pi\int_\alpha^\beta\sqrt{\det H(t)}\DD t
	.
\end{equation}
If $H(t)$ is invertible on some set of positive measure, this formula gives good information about the distribution of the
eigenvalues. On the other hand, if $\det H(t)=0$ almost everywhere, then it does not say anything other than that
$\sigma(A)$ is sparse compared to the integers. 

Denote $\xi_\alpha\DE\binom{\cos\alpha}{\sin\alpha}$. A Hamiltonian with zero determinant can always be written in the form 
\[
	H(t)=\Tr H(t)\cdot\xi_{\phi(t)}\xi_{\phi(t)}^T
	,
\]
where $\phi\DF I\to\bb R$ is a measurable function (determined up to integer multiples of $\pi$). We shall refer to $\phi$ as
the \emph{rotation angle} of $H$. 

The a basic question is how the distribution (density, asymptotics, etc.) of eigenvalues of 
selfadjoint realisations of \cref{Q8} relate to the rotation angle of $H$. 

Let us view this question from another angle which allows us to invoke function theory. 
We denote by $W(t,z)\DF I\times\bb C\to\bb C^{2\times 2}$ the unique solution of the initial value problem 
\[
	\left\{
	\begin{array}{l}
		\frac{\partial}{\partial t}W(t,z)J=zW(t,z)H(t),\quad t\in I\text{ a.e.}
		\\[2mm]
		W(\alpha,z)=I
	\end{array}
	\right.
\]
and call $W(t,z)$ the \emph{fundamental solution} of the system (for technical reasons we have passed to transposes, so
that the rows of $W(t,z)$ give the solutions of \cref{Q8}). The matrix $W_H(z)\DE W(\beta,z)$ is called the 
\emph{monodromy matrix} of the system. It is an entire function in the spectral parameter $z$. 

Understanding spectral properties amounts to understanding the monodromy matrix as an entire function 
because of the following central connection: 
there exists a selfadjoint realisation, call it $A_H$, such that all eigenvalues of $A_H$ are simple and 
\[
	\sigma(A_H)=\big\{x\in\bb R\DS w_{22}(x)=0\big\}\quad\text{where}\quad w_{22}(z)\DE(0,1)W_H(z)\binom 01
	.
\]
Hence, we have the immediate connection 
\begin{center}
	\emph{spectral distribution of selfadjoint realisations}
	\\
	$\longleftrightarrow$
	\\
	\emph{growth of the monodromy matrix}
\end{center}
and the correlation is that the slower the monodromy matrix grows, the less dense the spectrum will be. 

In the present paper we prove results which provide bounds for the growth of $W_H(z)$. 
Our main results are the three theorems described below. 
\\[1mm]
\emph{\Cref{Q6}:}\ 
In this theorem we provide a generic method to obtain upper bounds for $\log\|W_H(z)\|$. It should be seen as an improvement of
\cite[Theorem~1]{romanov:2017}. 
Formulation and proof are fairly similar, still \Cref{Q6} turns out to be a significant improvement of Romanov's Theorem. 
This can be witnessed even on the rough scale of exponential order, cf.\ \Cref{Q15}. 
As in Romanov's Theorem there is a lot of freedom when applying the result, and using this freedom in a clever
way is essential to obtain strong estimates. 
\\[1mm]
\emph{\Cref{Q14}:}\ 
We give an upper bound for the growth of $\log\|W_H\|$ for a Hamiltonian with continuous rotation angle. 
This is a perfect example for a (not too complicated) application of \Cref{Q6}. 
\\[1mm]
\emph{\Cref{Q49}:}\ 
In our third theorem we prove that the bound given in \Cref{Q14} is nearly sharp: we construct examples where the 
bound coming from \Cref{Q14} is equal to the maximum modulus up to a logarithmic factor. 
The proof requires major effort; among other things it relies on an auxiliary operator theoretic result which is of 
interest on its own right, cf.\ \Cref{Q25}. 

\smallskip
The sharpness result \Cref{Q49} is related to the following -- still open -- problem: Is it always possible to obtain the exact
growth of $W_H$ by an application of the bound obtained from Romanov's Theorem 
(naturally, in the form of the present improvement \Cref{Q6}) ?
There are several hints which indicate that the answer may be affirmative: \cite[Theorem~2]{romanov:2017} which
deals with diagonal Hamiltonians, \cite[Theorem~2.22]{pruckner.romanov.woracek:jaco} which deals with piecewise constant
Hamiltonians, and the present \Cref{Q49} which deals with continuous Hamiltonians. 

To close this introduction let us briefly describe the organisation of the content. We start with two 
sections containing auxiliary results. Those are needed only in the proof of the sharpness theorem, and therefore
the reader may skip Sections~2 and 3 until reaching \Cref{Q49}. Then we proceed to the stated main results:
in Section~4 we give the improvement of Romanov's Theorem, in Section~5 we apply it to obtain an upper bound for 
continuous rotation angles, and in Section~6 we prove sharpness for this case.

\section{Revisiting a lower bound for a Hamburger Hamiltonian}

In this section we discuss a lower bound for the growth of the monodromy matrix of Hamiltonians of a special form. 
This bound was used previously in \cite{pruckner.romanov.woracek:jaco,pruckner.woracek:srt} and (in a different language and
with a different proof) in \cite{berg.szwarc:2014}. A weaker variant appears already in \cite{livshits:1939}.

Recall that a \emph{Hamburger Hamiltonian} is a Hamiltonian of the form $H(t)=\Tr H(t)\cdot\xi_{\phi(t)}\xi_{\phi(t)}^T$ 
whose rotation angle $\phi(t)$ is piecewise constant with constancy intervals accumulating only at the right endpoint. 
More precisely:

\begin{Definition}
\label{Q22}
	Let $(l_j)_{j=1}^\infty$ be a summable sequence of positive numbers and $(\phi_j)_{j=1}^\infty$ be a sequence of real
	numbers. Set $L\DE\sum_{j=1}^\infty l_j$, and define a Hamiltonian $H_{l,\phi}$ on the interval $[0,L]$ as 
	\[
		H_{l,\phi}(t)\DE\xi_{\phi_j}\xi_{\phi_j}^T\quad\text{for }j\in\bb N\text{ and }
		\sum_{i=1}^{j-1}l_i\leq t<\sum_{i=1}^j l_i
		.
	\]
	A Hamiltonian $H_{l,\phi}$ thus can be pictured as 
	\begin{center}
	\begin{tikzpicture}[x=1.2pt,y=1.2pt,scale=0.8,font=\fontsize{8}{8}]
		\draw[thick] (10,30)--(215,30);
		\draw[dotted, thick] (215,30)--(270,30);
		\draw[thick] (10,25)--(10,35);
		\draw[thick] (70,25)--(70,35);
		\draw[thick] (120,25)--(120,35);
		\draw[thick] (160,25)--(160,35);
		\draw[thick] (190,25)--(190,35);
		\draw[thick] (210,25)--(210,35);
		\draw[thick] (270,25)--(270,35);
		\draw (40,44) node {${\displaystyle \xi_{\phi_1}\xi_{\phi_1}^T}$};
		\draw (95,44) node {${\displaystyle \xi_{\phi_2}\xi_{\phi_2}^T}$};
		\draw (140,44) node {${\displaystyle \xi_{\phi_3}\xi_{\phi_3}^T}$};
		\draw (177,43) node {${\cdots}$};
		\draw (-20,30) node {\large $H_{l,\phi}\!:$};
		\draw (10,15) node {${\displaystyle 0}$};
		\draw (70,15) node {${l_1}$};
		\draw (120,15) node {${l_1\!\!+\!l_2}$};
		\draw (160,15) node {${l_1\!\!+\!l_2\!\!+\!l_3}$};
		\draw (195,15) node {${\cdots}$};
		\draw (272,12) node {${\displaystyle L=\sum_{j=1}^\infty l_j}$};
	\end{tikzpicture}
	\end{center}
	We refer to the numbers $l_j$ and $\phi_j$ defining a Hamburger Hamiltonian as its \emph{lengths} and 
	\emph{angles}\footnote{Angles are determined only up to integer multiples of $\pi$}.
\end{Definition}

\noindent
This terminology is motivated from the connection with the Hamburger moment problem, see e.g.\ \cite{kac:1999}. 

The intuition concerning the growth of the monodromy matrix is that it grows slow if lengths decay fast, 
jumps of angles are small, and angles converge quickly. 
This reflects in the following result, which is the announced lower bound (it will also perfectly reflect in our later upper
bounds). 

\begin{Proposition}
\label{Q46}
	Let $H$ be a Hamburger Hamiltonian with lengths $(l_j)_{j=1}^\infty$ and angles $(\phi_j)_{j=1}^\infty$, and 
	assume that $\phi_1\not\equiv\frac\pi 2\!\!\!\mod\pi$. Set 
	\begin{equation}
	\label{Q47}
		F(z)\DE\sum_{n=0}^\infty \bigg[\prod_{j=1}^nl_{j+1}l_j\sin^2(\phi_{j+1}-\phi_j)\bigg]z^n
		,
	\end{equation}
	then 
	\[
		\log\Big(\max_{|z|=r}\|W_H(z)\|\Big)\geq\frac 12\log F(r^2)+\BigO(\log r)
		.
	\]
\end{Proposition}

\noindent
The assumption on $\phi_1$ is no loss of generality since adding a certain offset to the sequences of angles does not change the
function \cref{Q47} and changes $\log\|W_H(z)\|$ only up to a summand which is a $\BigO(\log|z|)$. 

The proof of \Cref{Q46} is obtained by repeating the ``Alternative proof of Proposition 2.15'' given in 
the extended preprint \cite[p.15]{pruckner.romanov.woracek:jaco-ASC}. 

\begin{proof}[Proof of \Cref{Q46} {\rm (cf.\ \cite{pruckner.romanov.woracek:jaco-ASC})}]
	For $t\geq 0,\phi\in\bb R$ and $z\in\bb C$ set 
	\[
		W_\phi(t,z)=I-zt\xi_\phi\xi_\phi^TJ
		,
	\]
	and note that $W_\phi(t,z)\xi_\phi=\xi_\phi$.
	Set $t_n\DE\sum_{j=1}^nl_j$, then the fundamental solution of $H$ is given as 
	\begin{multline*}
		W_H(t,z)=W_{\phi_1}(l_1,z)W_{\phi_2}(l_2,z)\ldots W_{\phi_{n-1}}(l_{n-1},z)W_{\phi_n}(t-t_{n-1},z)
		,
		\\
		\text{for}\quad n\in\bb N,\ t_{n-1}\leq t\leq t_n
		.
	\end{multline*}
	The function $(1,0)W_H(t,z)\xi_{\phi_n}$ is constant on the interval $[t_{n-1},t_n]$, and hence we can compute 
	(writing $W_H(z)=(w_{ij}(z))_{i,j=1}^2$) 
	\begin{align*}
		\frac{w_{12}(z)\ov{w_{11}(z)}-w_{11}(z)\ov{w_{12}(z)}}{z-\ov z}= &\, 
		\binom 10^T\int_0^L W_H(t,z)H(t)W_H(t,z)^*\DD t\binom 10
		\\
		& \mkern-85mu=\sum_{n=1}^\infty\int_{t_{n-1}}^{t_n}
		(1,0)W_H(t,z)\xi_{\phi_n}\cdot\xi_{\phi_n}^*W_H(t,z)^*\binom 10\DD t
		\\
		& \mkern-85mu=\sum_{n=1}^\infty \big|(1,0)W_H(t_{n-1},z)\xi_{\phi_n}\big|^2\cdot l_n
		.
	\end{align*}
	The function $p_n(z)\DE(1,0)W_H(t_{n-1},z)\xi_{\phi_n}$ is a polynomial of degree $n-1$ with real coefficients and 
	has only real zeroes. Therefore we have the estimate 
	\[
		|p_n(iy)|\geq y^{n-1}|c_{n-1}|,\quad y>0
		,
	\]
	where $c_{n-1}$ denotes the leading coefficient of $p_n(z)$. This coefficient computes as 
	\begin{align*}
		c_{n-1}\DE &\, (-1)^{n-1}l_1\cdot\ldots\cdot l_{n-1}\cdot(1,0)\cdot 
		\xi_{\phi_1}\xi_{\phi_1}^TJ\cdot\ldots\cdot\xi_{\phi_{n-1}}\xi_{\phi_{n-1}}^TJ\cdot\xi_{\phi_n}
		\\
		= &\, \Big(\prod_{j=1}^{n-1}l_j\Big)\cdot\cos\phi_1\cdot\Big(\prod_{j=1}^{n-1}\sin(\phi_{j+1}-\phi_j)\Big)
		.
	\end{align*}
	It follows that 
	\begin{align*}
		& |w_{11}(iy)|^2=
		\Big(\frac 1y\Im\frac{w_{12}(iy)}{w_{11}(iy)}\Big)^{-1}\cdot
		\frac{w_{12}(iy)\ov{w_{11}(iy)}-w_{11}(iy)\ov{w_{12}(iy)}}{2iy}
		\\
		&\mkern12mu \geq
		\Big(\frac 1y\Im\frac{w_{12}(iy)}{w_{11}(iy)}\Big)^{-1}\cdot
		\sum_{n=1}^\infty y^{2(n-1)}\bigg[\cos\phi_1\prod_{j=1}^{n-1}l_j\sin(\phi_{j+1}-\phi_j)\bigg]^2l_n
		\\
		&\mkern12mu = 
		\Big(\frac 1y\Im\frac{w_{12}(iy)}{w_{11}(iy)}\Big)^{-1}\cdot l_1\cos^2\phi_1\cdot
		\,\sum_{n=0}^\infty y^{2n}\bigg[\prod_{j=1}^nl_{j+1}l_j\sin^2(\phi_{j+1}-\phi_j)\bigg]
		.
	\end{align*}
	Each quotient of the entries of a line or a column of $W_H(z)$ is (up to a sign) a Herglotz function. We obtain 
	\begin{align*}
		\log\Big(\max_{|z|=r}\|W_H(z)\|\Big)\geq &\, \log\|W_H(ir)\|=\log|w_{11}(ir)|+\BigO(\log r)
		\\
		\geq &\, \frac 12\big(\log F(r^2)+\BigO(\log r)\big)+\BigO(\log r)
		.
	\end{align*}
\end{proof}

\noindent
In \Cref{Q76} we use this lower bound for Hamburger Hamiltonians whose lengths and angles are nicely behaving in the sense 
of regular variation (in Karamata's sense). For the theory of regular variation we refer to the monograph 
\cite{bingham.goldie.teugels:1989}; precise references will be given in course of the presentation. 
One can think of regularly varying functions as functions which behave roughly like a power. 
In this place, let us just recall the definition: a function $\ms f\DF[r_0,\infty)\to(0,\infty)$ defined on 
some ray is called \emph{regularly varying}, if there exists $\rho\in\bb R$ such that 
\[
	\forall\lambda>0\DP \lim_{r\to\infty}\frac{\ms f(\lambda r)}{\ms f(r)}=\lambda^\rho
	.
\]
The number $\rho$ is called the \emph{index} of $\ms f$, and we shall write $\Ind\ms f$ for it. 

One example which illustrates that regularly varying functions behave like powers in many respects is that they 
satisfy a variant of Stirlings approximation formula. We do not know an explicit reference and hence provide a proof.%
\footnote{
	Here, and throughout the paper, we shall use the following shorthand notations:
	\begin{align*}
		& f\sim g \DI \frac fg\to 1,\quad f\ll g \DI \frac fg\to 0
		\\
		& f\lesssim g \DI \exists C>0\DP f\leq Cg,\quad f\asymp g \DI \big(f\lesssim g\wedge g\lesssim f\big)
	\end{align*}
	}

\begin{lemma}
\label{Q44}
	Let $\ms f$ be regularly varying with index $\rho\in\bb R$. Then 
	\[
		\Big(\prod_{j=1}^n\ms f(j)\Big)^{\frac 1n}\sim\frac{\ms f(n)}{e^\rho}
		.
	\]
\end{lemma}
\begin{proof}
	Write $\ms f(r)=r^\rho\cdot\ms l(r)$ with $\ms l$ slowly varying. By Stirlings formula we have 
	\[
		\Big(\prod_{j=1}^nj^\rho\Big)^{\frac 1n}\sim \Big(\frac ne\Big)^\rho
		,
	\]
	hence we only have to deal with the slowly varying part. 

	By the representation theorem \cite[Theorem~1.3.1]{bingham.goldie.teugels:1989} we can write $\ms l$ as 
	\[
		\ms l(r)=c(r)\exp\Big(\int_1^r\frac{\epsilon(u)}u\DD u\Big)
	\]
	where $c$ and $\epsilon$ are bounded measurable functions such that $\lim_{r\to\infty}c(r)$ exists in $(0,\infty)$ and 
	$\lim_{r\to\infty}\epsilon(r)=0$. We obtain 
	\begin{equation}
	\label{Q45}
		\frac 1{\ms l(n)}\Big(\prod_{j=1}^n\ms l(j)\Big)^{\frac 1n}=
		\frac 1{c(n)}\Big(\prod_{j=1}^nc(j)\Big)^{\frac 1n}
		\cdot\,\exp\bigg(
		\frac 1n\sum_{j=1}^n\int\limits_1^j\frac{\epsilon(u)}u\DD u-\int\limits_1^n\frac{\epsilon(u)}u\DD u
		\bigg)
		.
	\end{equation}
	The first factor on the right side tends to $1$ because $c(r)$ has a positive and finite limit. 
	We estimate, for $j_0\geq 1$ and $n>j_0$,
	\begin{align*}
		\bigg|\frac 1n\sum_{j=1}^n\int_1^j\frac{\epsilon(u)}u &\, \DD u-\int_1^n\frac{\epsilon(u)}u\DD u\bigg|
		=\bigg|\frac 1n\sum_{j=1}^n\int_j^n\frac{\epsilon(u)}u\DD u\bigg|
		\\
		\leq &\, 
		\bigg|\frac 1n\sum_{j=1}^{j_0}\int_j^n\frac{\epsilon(u)}u\DD u\bigg|+
		\bigg|\frac 1n\sum_{j=j_0+1}^n\int_j^n\frac{\epsilon(u)}u\DD u\bigg|
		\\
		\leq &\, \frac{j_0\log n}n\cdot\sup_{r\in[1,\infty)}|\epsilon(r)|
		+\frac 1n\sum_{j=2}^n\int_j^n\frac 1u\DD u\cdot\sup_{r\in[j_0+1,\infty)}|\epsilon(r)|
		,
	\end{align*}
	and 
	\[
		\sum_{j=2}^n\int_j^n\frac 1u\DD u\leq\int\limits_1^n\int\limits_j^n\frac 1u\DD u\DD j
		=\int\limits_1^n\int\limits_1^u\frac 1u\DD j\DD u\leq n
		.
	\]
	Hence, also the second factor on the right side of \cref{Q45} tends to $1$. 
\end{proof}

\noindent
Further, recall an elementary lim-inf variant of the classical formula \cite[Theorem~I.2$^\prime$]{levin:1980} 
for the type w.r.t.\ a proximate order.

\begin{lemma}
\label{Q59}
	Let $A(z)=\sum_{n=0}^\infty a_nz^n$ be an entire function, let $r_0,s_0>0$ and 
	$\ms g\DF[r_0,\infty)\to[s_0,\infty)$ be an increasing bijection. Then 
	\[
		\liminf_{r\to\infty}\frac 1{\ms g(r)}\Big(\log\max\limits_{|z|=r}|A(z)|\Big)\geq
		\liminf_{n\to\infty}\log\Big(\ms g^{-1}(n)|a_n|^{\frac 1n}\Big)
		.
	\]
\end{lemma}
\begin{proof}
	For all $r>0$ and $n\in\bb N$ it holds that $\max_{|z|=r}|A(z)|\geq r^n|a_n|$, and in turn 
	\[
		\frac 1n\log\Big(\max_{|z|=r}|A(z)|\Big)\geq\log\big(r|a_n|^{\frac 1n}\big)
		.
	\]
	Using this for $r_n\DE\ms g^{-1}(n)$ gives 
	\[
		\liminf_{n\to\infty}\frac 1n\log\Big(\max_{|z|=r_n}|A(z)|\Big)\geq
		\liminf_{n\to\infty}\log\big(g^{-1}(n)|a_n|^{\frac 1n}\big)
	\]
	Let $r\geq r_1$ and take $n\in\bb N$ such that $r_n\leq r<r_{n+1}$. Then 
	\begin{align*}
		\frac 1{\ms g(r)}\log\Big(\max_{|z|=r}|A(z)|\Big)\geq &\, \frac 1{n+1}\log\Big(\max_{|z|=r_n}|A(z)|\Big)
		\\
		= &\, \frac n{n+1}\cdot\frac 1n\log\Big(\max_{|z|=r_n}|A(z)|\Big)
		,
	\end{align*}
	and it follows that 
	\[
		\liminf_{r\to\infty}\frac 1{\ms g(r)}\log\Big(\max_{|z|=r}|A(z)|\Big)\geq
		\liminf_{n\to\infty}\frac 1n\log\Big(\max_{|z|=r_n}|A(z)|\Big)
		.
	\]
\end{proof}

\noindent
Combining the above results yields the following lower bound for the maximum modulus of the monodromy matrix when lengths and
angles (in common) cannot have excessive downward drops. 

\begin{corollary}
\label{Q48}
	Let $H$ be a Hamburger Hamiltonian with lengths $(l_j)_{j=1}^\infty$ and angles $(\phi_j)_{j=1}^\infty$. 
	Let $\ms f$ be a regularly varying function, and choose $\ms g$ regularly varying with 
	$(\ms f\circ\ms g)(x)\sim(\ms g\circ\ms f)(x)\sim x$, see \cite[Theorem~1.5.12]{bingham.goldie.teugels:1989}. If 
	\begin{equation}
	\label{Q69}
		l_{j+1}l_j\sin^2(\phi_{j+1}-\phi_j)\gtrsim\frac 1{\ms f(j)}
		,\quad j\in\bb N
		,
	\end{equation}
	then 
	\begin{equation}
	\label{Q84}
		\log\max_{|z|=r}\|W_H(z)\|\gtrsim\ms g(r^2)
		.
	\end{equation}
\end{corollary}
\begin{proof}
	Since the sequence $(l_j)_{j=1}^\infty$ is summable, we also have 
	$\sum_{j=1}^\infty\big(\frac 1{\ms f(j)}\big)^{\frac 12}<\infty$. This implies that the index of $\ms f$, call it
	$\rho$, is at least $2$. 

	Passing from $\ms f$ to another regularly varying function $\tilde{\ms f}$ with $\ms f\asymp\tilde{\ms f}$ changes 
	$\ms g$ only up to ``$\asymp$'', and hence does not change the truth value of either \cref{Q69} or \cref{Q84}.
	We may use this freedom to assume without loss of generality that 
	\begin{Enumerate}
	\item $\ms f$ is an increasing bijection of $[1,\infty)$ onto itself,
	\item $\ms g=\ms f^{-1}$, 
	\item the assumption \cref{Q69} holds with ``$\geq$'' instead of ``$\gtrsim$''. 
	\end{Enumerate}
	Using \Cref{Q44} we obtain 
	\[
		\bigg(\prod_{j=1}^nl_{j+1}l_j\sin^2(\phi_{j+1}-\phi_j)\bigg)^{\frac 1n}
		\geq\bigg(\prod_{j=1}^n\frac 1{\ms f(j)}\bigg)^{\frac 1n}
		\sim\frac{e^\rho}{\ms f(n)}
		,
	\]
	and hence 
	\[
		\liminf_{n\to\infty}\log
		\bigg[\ms f(n)\cdot\bigg(\prod_{j=1}^nl_{j+1}l_j\sin^2(\phi_{j+1}-\phi_j)\big)^{\frac 1n}\bigg]
		\geq\rho\geq 2
		.
	\]
	It follows from \Cref{Q59} that the function $F(z)$ from \Cref{Q46} satisfies 
	\[
		\liminf_{r\to\infty}\frac 1{\ms g(r)}\log F(r)\geq 2
		.
	\]
	Note here that $F(z)$ has positive coefficients, and hence $\max_{|z|=r}|F(z)|=F(r)$. Now \Cref{Q46} gives 
	\[
		\log\Big(\max_{|z|=r}\|W_H(z)\|\Big)\geq\frac 12\log F(r^2)+\BigO(\log r)
		\geq\frac 12\ms g(r^2)+\BigO(\log r)\gtrsim\ms g(r^2)
		.
	\]
\end{proof}

\section{An auxiliary theorem from operator theory}

In this section we provide an auxiliary theorem about the operator model of a canonical system. 
It establishes a very intuitive fact, namely, that cutting out pieces of a Hamiltonian cannot increase the 
growth of the monodromy matrix. 

\subsection{The operator model of a canonical system}

To start with we briefly recall the definition and some properties of the operator model of the equation \cref{Q8}. 
Our standard reference in this respect is \cite{hassi.snoo.winkler:2000} and \cite[Chapter~7]{behrndt.hassi.snoo:2020}. 
The operator theory behind \cref{Q8} goes back to B.C.Orcutt \cite{orcutt:1969} and I.S.Kac \cite{kac:1984,kac:1986a} (see also 
\cite{kac:2002}), and in a different language to L.de~Branges \cite{debranges:1968}. Further recent references are 
\cite{remling:2018,romanov:1408.6022v1}.

Intervals where $H$ has constant nontrivial kernel require particular attention.

\begin{definition}
\label{Q77}
	Let $\phi\in\bb R$. A nonempty interval $(a,b)\subseteq I$ is called \emph{$H$-indivisible of type $\phi$}, if 
	\[
		H(t)=\Tr H(t)\cdot \xi_\phi\xi_\phi^T,\quad t\in(a,b)\text{ a.e.}
	\]
\end{definition}

\noindent
The type $\phi$ of an $H$-indivisible interval is unique up to integer multiples of $\pi$.
We shall assume throughout this section that the whole interval $I$ is not $H$-indivisible. This case is in some respects
trivial: the monodromy matrix is a linear polynomial. 

We denote by $L^2(H(t)\DD t)$ the usual $L^2$-space of equivalence classes of $2$-vector functions generated by the 
$2\!\times\! 2$-matrix measure $H(t)\DD t$, see e.g.\ \cite[p.1337--1346]{dunford.schwartz:1963}. 
To simplify notation, we shall always suppress explicit distinction between equivalence classes and their representants. 
However, one must keep in mind that sometimes it is important to make this distinction (for example when talking about boundary
values further below).

Now we can define the model space associated with a Hamiltonian $H$. 

\begin{definition}
\label{Q78}
	The \emph{model space $L^2(H)$} is the linear subspace of $L^2(H(t)\DD t)$ which consists of all functions 
	$f$ having the following property:
	\begin{Itemize}
	\item If $(a,b)$ is $H$-indivisible of type $\phi$, then $\xi_\phi^Tf(t)$ is constant a.e.\ on $(a,b)$.
	\end{Itemize}
\end{definition}

\noindent
The space $L^2(H)$ is a closed subspace of $L^2(H(t)\DD t)$, hence itself a Hilbert space, see e.g.\ 
\cite[Lemma~3.7]{hassi.snoo.winkler:2000}\footnote{%
	Caution: the notation in \cite{hassi.snoo.winkler:2000} is different. The space $L^2(H(t)\DD t)$ is what is 
	there called $L^2(H,\bb R^+)$, and our space $L^2(H)$ there is $L^2_s(H,\bb R^+)$.}. 

Next we define the minimal- and maximal- model operators. 

\begin{definition}
\label{Q79}
	Write $I=(\alpha,\beta)$. 
	The \emph{maximal-} and the \emph{minimal operators} $T_{\max}(H)$ and $T_{\min}(H)$ 
	are defined in terms of their graphs as 
	\begin{align*}
		T_{\max}(H)\DE &\, \bigg\{
		(f,g)\in L^2(H)\times L^2(H)\DS\ 
		\parbox{46mm}{\small $f$ has an absolutely continuous\\ representant with $f'=JHg$ a.e.}
		\bigg\}
		,
		\\
		T_{\min}(H)\DE &\, \bigg\{(f,g)\in T_{\max}(H)\DS\ 
		\parbox{50mm}{\small $f$ has an absolutely continuous\\ representant with $f(\alpha)=f(\beta)=0$}
		\bigg\}
		.
	\end{align*}
\end{definition}

\noindent
They have the following properties. 
\begin{Itemize}
\item For each $(f,g)\in T_{\max}(H)$ the first component has a unique absolutely continuous representant with $f'=JHg$. 
	Thus the boundary values $f(\alpha)$ and $f(\beta)$ are well-defined. 
\item An abstract Green's identity holds:
	\begin{multline*}
		\forall (f_1,g_1),(f_2,g_2)\in T_{\max}(H)\DP
		\\
		(g_1,f_2)_{L^2(H)}-(f_1,g_2)_{L^2(H)}=f_2(\alpha)^*Jf_1(\alpha)-f_2(\beta)^*Jf_1(\beta)
		.
	\end{multline*}
\item In some situations $T_{\max}(H)$ may be a multivalued operator. Despite this technical difficulty, it always 
	holds that $T_{\max}(H)=T_{\min}(H)^*$. 
\item $T_{\min}(H)$ is a closed symmetric operator, is completely nonselfadjoint, and has deficiency index $(2,2)$. 
\end{Itemize}
As a consequence of the above, selfadjoint extensions of $T_{\min}(H)$ can be described by boundary conditions at 
the left and right endpoints. We use the following two extensions: 
\begin{align*}
	B_H \DE &\, \big\{(f,g)\in T_{\max}(H)\DS f(\alpha)=0\big\},\quad R_H\DE B_H^{-1}
	,
	\\
	A_H \DE &\, \big\{(f,g)\in T_{\max}(H)\DS (1,0)f(\alpha)=(0,1)f(\beta)=0\big\}
	.
\end{align*}
The operator $R_H$ is the Volterra integral operator
\[
	(R_Hf)(t)\DE\int_\alpha^tJH(t)f(t)\DD t
	,\quad f\in L^2(H)
	,
\]
while $A_H$ is selfadjoint. Note that $B_H$ and $A_H$ are invertible since 
$\ker T_{\max}(H)=\Span\big\{\binom 10,\binom 01\big\}$, and that $A_H^{-1}$ is a rank-one perturbation of $R_H$. 

\subsection{Cutting out pieces of a Hamiltonian}

The theorem announced at the beginning of this section reads as follows.

\begin{theorem}
\label{Q25}
	Let $H$ be a Hamiltonian on $I=[\alpha,\beta]$. Let $\Delta\subseteq[\alpha,\beta]$ be a (Lebesgue-) measurable set
	with positive measure, 
	and assume that for every $H$-indivisible interval $(a,b)\subseteq[\alpha,\beta]$ either $(a,b)\cap\Delta$ or 
	$(a,b)\setminus\Delta$ has measure zero. Set 
	\begin{align*}
		\lambda(t)\DE &\, 
		\int_\alpha^t\mathds{1}_{\Delta}(u)\DD u,\ t\in[\alpha,\beta],\qquad \tilde L\DE\lambda(\beta)
		,
		\\
		\kappa(s)\DE &\, 
		\min\big\{t\in[\alpha,\beta]\DS \lambda(t)=s\big\},\quad s\in[0,\tilde L]
		,
		\\
		\tilde H\DE &\, 
		H\circ\kappa
		.
	\end{align*}
	Then the following statements hold. 
	\begin{Enumerate}
	\item $\tilde H$ is a Hamiltonian on $[0,\tilde L]$, and satisfies 
		\begin{equation}
		\label{Q43}
			(\tilde H\circ\lambda)\cdot\mathds{1}_{\Delta}=H\cdot\mathds{1}_{\Delta}
			\text{ a.e.}
		\end{equation}
	\item The map $V$ acting as 
		\[
			V\DF f\mapsto (f\circ\lambda)\cdot\mathds{1}_{\Delta}
		\]
		induces an isometry of $L^2(\tilde H)$ into $L^2(H)$. 
	\item Denote by $M_{\mathds{1}_{\Delta}}$ the multiplication operator with $\mathds{1}_{\Delta}$. Then we have that 
		$\Ran M_{\mathds{1}_{\Delta}}R_HV\subseteq\Ran V$ and $R_{\tilde H}=V^{-1}M_{\mathds{1}_{\Delta}}R_HV$. 
		\begin{center}
		\begin{tikzcd}[column sep=normal]
			L^2(H) \arrow{r}{R_H} & L^2(H) \arrow{r}{M_{\mathds{1}_{\Delta}}} 
			& \Ran(M_{\mathds{1}_{\Delta}}R_HV) \arrow{d}{\subseteq}
			\\
			L^2(\tilde H) \arrow{r}[swap]{R_{\tilde H}} \arrow{u}{V}
			& L^2(\tilde H) \arrow[r, bend left=15, "V"] \ar[r, phantom, "\cong"] 
			& \Ran V \ar[l, bend left=15, "V^{-1}"]
		\end{tikzcd}
		\end{center}
	\end{Enumerate}
\end{theorem}
\begin{proof}
	The proof of \Enumref{1} relies on some measure theoretic considerations. 
	Let us denote the maximal constancy
	intervals which contain more than one point (if any) as $\Delta_j$. There exist at most countably many such intervals
	and $\Delta_j\cap\Delta$ is a zero set for all $j$. We now show that (the complement is understood in $[\alpha,\beta]$)
	\[
		\lambda(\Delta^c)\text{ is a zero set}
		.
	\]
	Since $\lambda$ is absolutely continuous, the set $\lambda(\Delta^c)$ is (Lebesgue-) measurable. The change of
	variables formula gives 
	\[
		\int_0^{\tilde L}\mathds{1}_{\lambda(\Delta^c)}(s)\DD s
		=\int_\alpha^\beta\big(\mathds{1}_{\lambda(\Delta^c)}\circ\lambda\big)(t)\cdot\mathds{1}_{\Delta}(t)\DD t
		=\int_\alpha^\beta\mathds{1}_{\lambda^{-1}(\lambda(\Delta^c))\cap\Delta}(t)\DD t
		.
	\]
	We have 
	\[
		\lambda^{-1}(\lambda(\Delta^c))=\Delta^c\cup\bigcup\big\{\Delta_j\DS \Delta_j\cap\Delta\neq\emptyset\big\}
		,
	\]
	and hence the integral on the right vanishes. 

	In the second step we show that the function $\kappa\DF[0,\tilde L]\to[\alpha,\beta]$, which is defined as 
	\[
		\kappa(s)\DE\min\big\{t\in[\alpha,\beta]\DS \lambda(t)=s\big\},\quad s\in[0,\tilde L]
		,
	\]
	is Lebesgue-to-Lebesgue measurable. Clearly, $\kappa$ is nondecreasing and a right inverse of $\lambda$. 
	Monotonicity implies that it is Borel-to-Borel measurable. Let $E$ be a Lebesgue measurable subset of $[\alpha,\beta]$, 
	and choose Borel sets $A,B\subseteq[\alpha,\beta]$ with $A\subseteq E\subseteq B$ and $B\setminus A$ being a zero set. 
	Then $\kappa^{-1}(A)\subseteq\kappa^{-1}(E)\subseteq\kappa^{-1}(B)$ and 
	\[
		\kappa^{-1}(B)\setminus\kappa^{-1}(A)=\kappa^{-1}(B\setminus A)\subseteq\lambda(B\setminus A)
		.
	\]
	The set on the right is a zero set since $\lambda$ is absolutely continuous, and we conclude that 
	$\kappa^{-1}(E)$ is Lebesgue measurable. 

	Now we define 
	\[
		\tilde H\DE H\circ\kappa\DF[0,\tilde L]\to\bb R^{2\times 2}
		.
	\]
	Obviously, $\tilde H$ takes nonnegative matrices as values and is (Lebesgue-) measurable. 
	Moreover, we have 
	\[
		\big\{t\in[\alpha,\beta]\DS (\kappa\circ\lambda)(t)\neq t\big\}\subseteq\bigcup_j\Delta_j
		,
	\]
	and hence 
	\[
		(\tilde H\circ\lambda)(t)\mathds{1}_{\Delta}(t)=H(\kappa\circ\lambda(t))\mathds{1}_{\Delta}(t)
		=H(t)\mathds{1}_{\Delta}(t)\quad\text{a.e.}
	\]
	This is \cref{Q43}.
	We need to check that $\tilde H$ is a Hamiltonian. Let $B\subseteq[0,\tilde L]$ be measurable, then 
	\begin{align*}
		\int_0^{\tilde L}\Tr\tilde H(s)\mathds{1}_B(s)\DD s
		= &\, \int_\alpha^\beta\Tr(\tilde H\circ\lambda)(t)(\mathds{1}_B\circ\lambda)(t)\mathds{1}_{\Delta}(t)\DD t
		\\
		= &\, \int_\alpha^\beta\Tr H(t)\mathds{1}_{\lambda^{-1}(B)\cap\Delta}\DD t
		.
	\end{align*}
	Choosing $B=[0,\tilde L]$ already shows that $\tilde H$ is integrable. Assume now that $B$ is some set with positive
	measure. Since $\Tr H(t)>0$ a.e., measurability of the integrand in the last integral implies that the set 
	$\lambda^{-1}(B)\cap\Delta$ is measurable. Moreover, $B\setminus\lambda(\lambda^{-1}(B)\cap\Delta)$ is a zero set since
	it is contained in $\lambda(\Delta^c)$. Therefore $\lambda^{-1}(B)\cap\Delta$ must have positive measure, and the
	integral on the right is positive. 

	We come to the proof of \Enumref{2}.
	The first step is to observe that $V$ is isometric. This follows simply by making a change of variable. 
	Let $f\DF[0,\tilde L]\to\bb C^2$ be any measurable function, then we have
	\begin{align*}
		\int_0^{\tilde L} f(s)^*\tilde H(s)f(s)\DD s= &\, 
		\int_\alpha^\beta (f\circ\lambda)(t)^*(\tilde H\circ\lambda)(t)(f\circ\lambda)(t)\cdot\mathds{1}_{\Delta}(t)\DD t
		\\
		= &\, \int_\alpha^\beta (Vf)(t)^*H(t)(Vf)(t)\DD t
		.
	\end{align*}
	Note that isometry implies
	\begin{equation}
	\label{Q31}
		\tilde Hf_1=\tilde Hf_2\ \text{a.e.}
		\quad\Longrightarrow\quad
		H(Vf_1)=H(Vf_2)\ \text{a.e.}
	\end{equation}
	We have to check the constancy condition from \Cref{Q78} for indivisible intervals. 
	Let $f\DF[0,\tilde L]\to\bb C^2$ be a measurable function which satisfies the condition for $\tilde H$. 
	We have to show that $Vf$ satisfies it for $H$. 

	Let $(a,b)\subseteq[\alpha,\beta]$ be an $H$-indivisible interval, and let $\phi$ be its type. 
	By the assumption of the theorem, either $(a,b)\cap\Delta$ or $(a,b)\setminus\Delta$ is a zero set. In the first case,
	we have $(Vf)(t)=0$ for $t\in(a,b)$ a.e., and are done. Consider the second case. Then 
	$\mathds{1}_{\Delta}(t)=1$ for $t\in(a,b)$ a.e., and hence 
	\begin{equation}
	\label{Q80}
		\tilde H(\lambda(t))=H(t)=\Tr H(t)\cdot\xi_\phi\xi_\phi^T
		,\quad t\in(a,b)\ \text{a.e.}
	\end{equation}
	Since $\lambda$ is absolutely continuous and nondecreasing, we have $(\lambda(a),\lambda(b))\subseteq\lambda((a,b))$ and 
	the image of the exceptional set in \cref{Q80} is a zero set. Hence, 
	\[
		\tilde H(s)=\Tr\tilde H(s)\cdot\xi_\phi\xi_\phi^T
		,\quad s\in(\lambda(a),\lambda(b))\text{ a.e.}
	\]
	This means that $(\lambda(a),\lambda(b))$ is $\tilde H$-indivisible of type $\phi$, and hence that $\xi_\phi^Tf(s)$ is 
	constant on $(\lambda(a),\lambda(b))$ a.e. Say, we have $\xi_\phi^Tf(s)=\gamma$ for a.a.\ 
	$s\in(\lambda(a),\lambda(b))$. It follows that 
	$\tilde H(s)f(s)=\tilde H(s)(\gamma\xi_\phi)$ for a.a.\ $s\in(\lambda(a),\lambda(b))$, in other words, the functions 
	\[
		f_1\DE\mathds{1}_{(\lambda(a),\lambda(b))}f,\quad
		f_2\DE\mathds{1}_{(\lambda(a),\lambda(b))}(\gamma\xi_\phi)
	\]
	satisfy $\tilde Hf_1=\tilde Hf_2$ for a.a.\ $s\in[0,\tilde L]$. Applying \cref{Q31} yields
	\[
		H(t)\big[(\mathds{1}_{(\lambda(a),\lambda(b))}\circ\lambda)(t)(f\circ\lambda)(t)\mathds{1}_{\Delta}(t)\big]
		=H(t)\big[(\mathds{1}_{(\lambda(a),\lambda(b))}\circ\lambda)(t)\gamma\xi_\phi\mathds{1}_{\Delta}(t)\big]
	\]
	for a.a.\ $t\in[\alpha,\beta]$. Since $(a,b)\setminus\Delta$ is a zero set, we have $\mathds{1}_{\Delta}(t)=1$ for a.a.\ 
	$t\in(a,b)$. The function $\lambda$ is strictly increasing on $(a,b)$, and hence
	$(\mathds{1}_{(\lambda(a),\lambda(b))}\circ\lambda)(t)=1$ for all $t\in(a,b)$. It follows that 
	\[
		\Tr H(t)\cdot\xi_\phi^T(Vf)(t)=\Tr H(t)\cdot \xi_\phi^T(\gamma\xi_\phi)
		,\quad t\in(a,b)\text{ a.e.}
		,
	\]
	and hence $\xi_\phi^T(Vf)(t)=\gamma$ again for $t\in(a,b)$ a.e.
	
	Finally, we come to the proof of \Enumref{3}. First note that, by our assumption on indivisible intervals, 
	$M_{\mathds{1}_\Delta}$ maps $L^2(H)$ into itself (in fact, is an orthogonal projection). 
	Now let $f\in L^2(\tilde H)$. Then 
	\begin{align*}
		\big[(R_H\circ V)(f)\big](t)= &\, \int_\alpha^t JH(u)\cdot(f\circ\lambda)(u)\mathds{1}_{\Delta}(u)\DD u
		\\
		= &\, \int_\alpha^t J(\tilde H\circ\lambda)(u)(f\circ\lambda)(u)\mathds{1}_{\Delta}(u)\DD u
		\\
		= &\, \int_0^{\lambda(t)} J\tilde H(r)f(r)\DD r
		=\big[(R_{\tilde H}f)\circ\lambda\big](t)
		.
	\end{align*}
	We see that $M_{\mathds{1}_{\Delta}}R_HV=VR_{\tilde H}$, and the assertion follows.
\end{proof}

\noindent
Let us note that $\tilde H$ defined above is the unique Hamiltonian with \cref{Q43}. 
To see this, assume we have $\hat H$ with \cref{Q43}. Then 
\[
	\hat H\cdot(\mathds{1}_{\Delta}\circ\kappa)=\big[(\hat H\circ\lambda)\cdot\mathds{1}_{\Delta}\big]\circ\kappa
	=(H\cdot\mathds{1}_{\Delta})\circ\kappa=\tilde H\cdot(\mathds{1}_{\Delta}\circ\kappa)
	.
\]
We have $\mathds{1}\circ\kappa=\mathds{1}_{\kappa^{-1}(\Delta)}$, and since
$\kappa^{-1}(\Delta^c)\subseteq\lambda(\Delta^c)$ this is equal to $1$ a.e.

Passing to growth properties of $W_H$ can easily be done using the usual function theoretic tools. 

\begin{corollary}
\label{Q17}
	Consider the situation described in \Cref{Q25}. 
	Moreover, let $\ms f$ be a regularly varying function with index $\rho\in(0,1)$. 
	Then 
	\[
		\limsup_{|z|\to\infty}\frac{\log\|W_{\tilde H}(z)\|}{\ms f(|z|)}\lesssim
		\limsup_{|z|\to\infty}\frac{\log\|W_H(z)\|}{\ms f(|z|)}
		.
	\]
	The constant implicit in this relation depends only on $\rho$. 
\end{corollary}
\begin{proof}
	For a compact operator $T$ we denote by $s_n(T)$ its $n$-th s-number and let $n_T(r)$ be the counting function 
	\[
		n_T(r)\DE\#\big\{n\DS s_n(T)\geq\frac 1r\big\},\quad r>0
		.
	\]
	Due to \Cref{Q25} we have $s_n(R_{\tilde H})\leq s_n(R_H)$ for all $n$, 
	and hence $n_{R_{\tilde H}}(r)\leq n_{R_H}(r)$ for all $r>0$. 

	The operator $A_H^{-1}$ is a rank-one perturbation of $R_H$, and the same for $A_{\tilde H}^{-1}$ and 
	$R_{\tilde H}$. Hence, we have 
	\[
		n_{A_{\tilde H}^{-1}}(r)\leq n_{R_{\tilde H}}(r)+1\leq n_{R_H}(r)+1\leq n_{A_H^{-1}}(r)+2
		.
	\]
	The spectrum of $A_H$ coincides with the zero set of the entire function $w_{22}(z)\DE(0,1)W_H(z)\binom 01$, and 
	the spectrum of $A_{\tilde H}$ with the zero set of $\tilde w_{22}(z)\DE(0,1)\tilde W(z)\binom 01$. Thus we have 
	(now using the notation $n_f(r)$ for the counting function of the zeroes of an entire function $f$) 
	\[
		n_{\tilde w_{22}}(r)\leq n_{w_{22}}(r)+2,\quad r>0
		.
	\]
	Due to \cite[Proposition~7.4.1]{bingham.goldie.teugels:1989} we can assume w.l.o.g.\ that $\ms f$ is a proximate order.
	Now \cite[Theorem~I.17]{levin:1980} is applicable, and yields 
	\[
		\limsup_{|z|\to\infty}\frac{\log|\tilde w_{22}(z)|}{\ms f(|z|)}
		\asymp\limsup_{r\to\infty}\frac{n_{\tilde w_{22}}(r)}{\ms f(r)}
		\leq\limsup_{r\to\infty}\frac{n_{w_{22}}(r)}{\ms f(r)}
		\asymp\limsup_{|z|\to\infty}\frac{\log|\tilde w_{22}(z)|}{\ms f(|z|)}
		.
	\]
	By the proof of \cite[Theorem~I.17]{levin:1980}, 
	the constants implicit in this relation depend only on the index of $\ms f$. 
\end{proof}

\section{A generic estimate from above}

In the below theorem we provide a method to estimate the monodromy matrix of a canonical system. 
This result is an improvement of a theorem due to R.Romanov in \cite{romanov:2017}. 
The proof follows the very same idea as \cite[Theorem~1]{romanov:2017} and -- despite the result being stronger -- 
the argument is equally simple: it merely uses multiplicativity of the fundamental solution and 
Gr\"onwall's lemma. Similar as for \cite[Theorem~1]{romanov:2017} the power of \Cref{Q6} is its flexibility. Applying it in 
a clever way is at least as important as the theorem itself. 

For practical reasons we throughout use the spectral norm on $\bb C^{2\times 2}$. 
This norm has the advantage to be invariant under unitary transformations. 

\begin{Theorem}
\label{Q6}
	Let $H$ be a Hamiltonian on a compact interval $I$ with $\det H=0$ a.e., and write
	$H(t)=\Tr H(t)\cdot\xi_{\phi(t)}\xi_{\phi(t)}^T$ with a measurable function $\phi\DF I\to\bb R$. 
	Assume we are given
	\begin{Ilist}
	\item a partition $(y_0,\ldots,y_N)$ of $I$, i.e., 
		\[
			N\in\bb N,\quad \min I=y_0<y_1<\ldots<y_N=\max I
			,
		\]
	\item rotation parameters $\psi_1,\ldots,\psi_N\in\bb R$,
	\item distortion parameters $a_1,\ldots,a_N\in(0,1]$,
	\end{Ilist}
	and set 
	\begin{align*}
		A_1 \DE &\, 
		\sum_{j=1}^N a_j^2\int_{y_{j-1}}^{y_j}\cos^2\big(\phi(t)-\psi_j\big)\cdot\Tr H(t)\DD t
		,
		\\
		A_2 \DE &\, 
		\sum_{j=1}^N \frac 1{a_j^2}\int_{y_{j-1}}^{y_j}\sin^2\big(\phi(t)-\psi_j\big)\cdot\Tr H(t)\DD t
		,
		\\
		A_3 \DE &\, 
		\sum_{j=1}^{N-1} \log\bigg(
		\max\Big\{\frac{a_j}{a_{j+1}},\frac{a_{j+1}}{a_j}\Big\}\cdot\big|\cos\big(\psi_j-\psi_{j+1}\big)\big|
		+\frac{|\sin(\psi_j-\psi_{j+1})|}{a_ja_{j+1}}
		\bigg)
		,
		\\
		A_4 \DE &\, 
		-\log a_1-\log a_N
		.
	\end{align*}
	Then 
	\begin{equation}
	\label{Q7}
		\forall z\in\bb C\DP \log\|W_H(z)\|\leq|z|\cdot(A_1+A_2)+A_3+A_4
		,
	\end{equation}
	where $\|\Dummy\|$ denotes the spectral norm on $\bb C^{2\times 2}$.
\end{Theorem}

\noindent
The following remark is essential for successful application of \Cref{Q6}. 

\begin{Remark}
\label{Q70}
	On first sight the estimate \cref{Q7} may seem quite useless. We know from the Krein-de~Branges formula that 
	$W_H(z)$ is of minimal exponential type, and hence of course an estimate $\log\|W_H(z)\|\lesssim |z|$ holds. 
	The significance of \Cref{Q6} lies in a quantitative aspect. Namely, \cref{Q7} holds for \emph{all} choices of data 
	$y_j,\psi_j,a_j$ for \emph{all} complex numbers $z$. 

	Now reverse the viewpoint. Consider $z$ as fixed, choose the data $y_j(z),\psi_j(z),a_j(z)$ in dependence of $z$, 
	and use \cref{Q7} only for the given point $z$. If we manage to make the $z$-dependent choice of data in such a way that
	$A_1$ and $A_2$ decay when $|z|$ increases to $\infty$, and that $A_3$ and $A_4$ do not grow too fast, we may get a
	bound for $\log\|W_H(z)\|$ which is significantly smaller than $|z|$. 
\end{Remark}

\noindent
For the proof of \Cref{Q6} we start with an application of Gr\"onwall's lemma.

\begin{lemma}
\label{Q1}
	Let $H$ be a Hamiltonian on a compact interval $I$.
	Assume we are given a partition $(y_0,\ldots,y_N)$ of $I$ and matrices $\Omega_1,\ldots,\Omega_N\in\GL(2,\bb R)$. 
	Then (for any submultiplicative norm)
	\begin{equation}
	\label{Q2}
		\|W_H(z)\|\leq
		\exp\bigg(|z|\sum_{j=1}^N \int_{y_{j-1}}^{y_j}\|\Omega_j H(t)J\Omega_j^{-1}\|\DD t\bigg)
		\|\Omega_1^{-1}\|\|\Omega_N\|\prod_{j=1}^{N-1}\|\Omega_j\Omega_{j+1}^{-1}\|
		.
	\end{equation}
\end{lemma}
\begin{proof} 
	For each $j\in\{1,\ldots,N\}$ let $W_j(t,z)$ be the fundamental solution of $H|_{[y_{j-1},y_j]}$, 
	and let $W_j\DE W_j(y_j,z)$ be the corresponding monodromy matrices. Then 
	\[
		W_H(z)=W_1(z)\cdot W_2(z)\cdot\ldots\cdot W_N(z) 
		.
	\]
	We insert the matrices $\Omega_j$ and get
	\[
		W_H(z)=\Omega_1^{-1}\cdot\big(\Omega_1 W_1(z)\Omega_1^{-1}\big)\cdot\Omega_1\Omega_2^{-1}\cdot
		\ldots\cdot\big(\Omega_N W_N(z)\Omega_N^{-1}\big)\cdot\Omega_N
		.
	\]
	Applying Gr\"onwall's lemma to the differential equation
	\[
		\frac{\partial}{\partial x}\Omega_j W_j(x,z)\Omega_j^{-1}=
		-z\cdot\Omega_j W_j(x,z)\Omega_j^{-1}\cdot\Omega_jH(x)J\Omega_j^{-1}
		,\quad x\in[y_{j-1},y_j]
		,
	\]
	yields that
	\[
		\|\Omega_j W_j(z)\Omega_j^{-1}\|\leq
		\exp\Big(|z|\int_{y_{j-1}}^{y_j}\|\Omega_j H(t)J\Omega_j^{-1}\|\DD t\Big) 
		.
	\]
	The assertion of the lemma follows.
\end{proof}

\noindent
There happens no loss in precision when using only matrices $\Omega_j$ of a particular form.

\begin{Definition}
\label{Q4}
	For $a,b>0$ set $D(a,b)\DE(\begin{smallmatrix} a & 0 \\ 0 & b \end{smallmatrix})$ and denote, for $a>0$ and 
	$\psi\in\bb R$, 
	\[
		\Omega(a,\psi)\DE D(a,a^{-1})\exp(-\psi J)
		=\begin{pmatrix} a & 0 \\ 0 & a^{-1}\end{pmatrix}
		\begin{pmatrix} \cos\psi & -\sin\psi \\ \sin\psi & \cos\psi\end{pmatrix}
		.
	\]
\end{Definition}

\noindent
Geometrically, the matrix $\Omega(a,\psi)$ is a rotation followed by a distortion. 

\begin{Remark}
\label{Q3}
	To see that we may restrict to matrices of the form $\Omega(a,\psi)$,
	observe that the right side of \cref{Q2} remains unchanged when the matrices $\Omega_j$ are multiplied 
	with real nonzero scalars $\alpha_j$ or multiplied from the left with matrices $C_j\in\GL(2,\bb R)$ satisfying
	$\|C_j\|=\|C_j^{-1}\|=1$.
	Using these two transformations, every matrix $\Omega\in\GL(2,\bb R)$ can be brought to the form $\Omega(a,\psi)$. 
\end{Remark}

\noindent
In the next lemma we compute the relevant norms for matrices $\Omega(a,\psi)$. 

\begin{lemma}
\label{Q5}
	Let $a,b>0$ and $\psi,\phi\in\bb R$. 
	\begin{Enumerate}
	\item \Dis{%
		\|\Omega(a,\psi)\|=\|\Omega(a,\psi)^{-1}\|=\max\big\{a,a^{-1}\big\}
		.
		}
	\item \Dis{%
		\|\Omega(a,\psi)\xi_\phi\xi_\phi^TJ\Omega(a,\psi)^{-1}\|
		=a^2\cos^2(\phi-\psi)+\frac 1{a^2}\sin^2(\phi-\psi)
		.
		}
	\item Set 
		\[
			v_+\DE\binom{\max\{\frac ab,\frac ba\}|\cos(\phi-\psi)|}{\max\{ab,\frac 1{ab}\}|\sin(\phi-\psi)|}
			,
			v_-\DE\binom{\min\{\frac ab,\frac ba\}|\cos(\phi-\psi)|}{\min\{ab,\frac 1{ab}\}|\sin(\phi-\psi)|}
			,
		\]
		and denote by $\|\Dummy\|_p$, $p\in\{1,2\}$, the $p$-norm on $\bb R^2$. Then 
		\begin{align*}
			\|v_+\|_2^2\leq &\, \|\Omega(a,\psi)\Omega(b,\phi)^{-1}\|^2
			\\
			= &\, 1+\|v_+-v_-\|_2\cdot\frac{\|v_+-v_-\|_2+\|v_++v_-\|_2}2
			\\
			\leq &\, \|v_+\|_1^2\leq 2\|v_+\|_2^2
			.
		\end{align*}
	\end{Enumerate}
\end{lemma}
\begin{proof}
	For the proof of \Enumref{1} it is enough to note that $\exp(-\psi J)$ is unitary. This implies that 
	\[
		\|\Omega(a,\psi)\|=\|D(a,a^{-1})\exp(-\psi J)\|=\|D(a,a^{-1})\|=\max\{a,a^{-1}\}
		,
	\]
	and the analogous formula for $\Omega(a,\psi)^{-1}$. 

	We come to the proof of \Enumref{2}. Note the relations $J\exp(\phi J)=\exp(\phi J)J$ and 
	$\xi_\phi\xi_\phi^T=\exp(\phi J)(\begin{smallmatrix} 1 & 0 \\ 0 & 0 \end{smallmatrix})\exp(-\phi J)$, 
	which are easily verified. Moreover, set $\sigma\DE\phi-\psi$. Then 
	\begin{align*}
		B\DE &\, \Omega(a,\psi)\xi_{\phi}\xi_{\phi}^TJ\Omega(a,\psi)^{-1} 
		\\
		= &\, D(a,a^{-1})\exp(-\psi J)\xi_{\phi}\xi_{\phi}^TJ\exp(\psi J)D\big(a^{-1},a\big)
		\\
		= &\, D(a,a^{-1})\exp(-\psi J)\exp(\phi J)
		\big(\begin{smallmatrix} 1 & 0 \\ 0 & 0 \end{smallmatrix}\big)
		\exp(-\phi J)\exp(\psi J)JD\big(a^{-1},a\big)
		\\
		= &\, D(a,a^{-1})\xi_\sigma\xi_\sigma^TJD\big(a^{-1},a\big) 
		=\begin{pmatrix} 
			\cos(\sigma)\sin(\sigma) & -a^2\cos^2(\sigma) 
			\\
			\frac 1{a^2}\sin^2(\sigma) & -\cos(\sigma)\sin(\sigma)
		\end{pmatrix}
		.
	\end{align*}
	A direct computation shows
	\[
		B^TB=
		\begin{pmatrix} 
			\frac{1}{a^4}\sin^4(\sigma)+\cos^2(\sigma)\sin^2(\sigma) & \ast
			\\
			\ast & a^4\cos^4(\sigma)+\cos^2(\sigma)\sin^2(\sigma)  
		\end{pmatrix}
		,
	\]
	from which we see that $\Tr(B^T B)=(a^2\cos^2(\sigma)+\frac{1}{a^2}\sin^2(\sigma))^2$.
	Since $\det B=0$, we have $\|B\|=\sqrt{\Tr(B^TB)}$.

	Finally, we turn to \Enumref{3}. We compute 
	\begin{align*}
		C\DE &\, \Omega(a,\psi)\Omega(b,\phi)^{-1} 
		=D(a,a^{-1})\exp((\phi-\psi)J)D(b^{-1},b)
		\\[1mm]
		= &\, 
		\begin{pmatrix} 
			\frac ab\cos\sigma & -ab\sin\sigma
			\\
			\frac 1{ab}\sin\sigma & \frac ba\cos\sigma
		\end{pmatrix}
		=\cos\sigma\begin{pmatrix} \frac ab & 0 \\ 0 & \frac ba \end{pmatrix}
		+\sin\sigma\begin{pmatrix} 0 & -ab \\ \frac 1{ab} & 0 \end{pmatrix}
		.
	\end{align*}
	The asserted estimate from above follows:
	\begin{align*}
		\|C\|\leq &\, |\cos\sigma|\cdot\bigg\|\begin{pmatrix} \frac ab & 0 \\ 0 & \frac ba \end{pmatrix}\bigg\|
		+|\sin\sigma|\cdot\bigg\|\begin{pmatrix} 0 & -ab \\ \frac 1{ab} & 0 \end{pmatrix}\bigg\|
		\\
		= &\, |\cos\sigma|\cdot\max\Big\{\frac ab,\frac ba\Big\}+|\sin\sigma|\cdot\max\Big\{ab,\frac 1{ab}\Big\}
		=\|v_+\|_1
		.
	\end{align*}
	A calculation shows
	\[
		C^TC=
		\begin{pmatrix}
			\cos^2(\sigma)\big(\frac ab\big)^2+\sin^2(\sigma)\big(\frac 1{ab}\big)^2 & \ast
			\\
			\ast & \cos^2(\sigma)\big(\frac ba\big)^2+\sin^2(\sigma)(ab)^2
		\end{pmatrix}
		,
	\]
	and we see that 
	\[
		\Tr(C^TC)=
		\cos^2\sigma\cdot\Big[\big(\frac ab\big)^2+\big(\frac ba\big)^2\Big]
		+\sin^2\sigma\cdot\Big[(ab)^2+\big(\frac 1{ab}\big)^2\Big]
		.
	\]
	We have $\det(C^TC)=1$, and hence the eigenvalues of $C^TC$ are the solutions of the equation 
	\begin{equation}
	\label{Q10}
		\lambda+\frac 1\lambda=\Tr(C^TC)
		.
	\end{equation}
	To shorten notation, set $\tau\DE\Tr(C^TC)$. Computing the larger of the solutions of \cref{Q10} gives 
	\begin{align*}
		\|C^TC\|= &\, \frac 12\big(\tau+\sqrt{\tau^2-4}\big)
		=1+\frac 12\big((\tau-2)+\sqrt{(\tau-2)(\tau+2)}\big)
		\\
		= &\, 1+(\tau-2)^{\frac 12}\cdot\frac{(\tau-2)^{\frac 12}+(\tau+2)^{\frac 12}}{2}
		.
	\end{align*}
	Now note that 
	\begin{align*}
		\tau-2 = &\, \cos^2\sigma\cdot\Big(\frac ab-\frac ba\Big)^2+\sin^2\sigma\cdot\Big(ab-\frac 1{ab}\Big)^2
		=\|v_+-v_-\|_2^2
		,
		\\
		\tau+2 = &\, \cos^2\sigma\cdot\Big(\frac ab+\frac ba\Big)^2+\sin^2\sigma\cdot\Big(ab+\frac 1{ab}\Big)^2
		=\|v_++v_-\|_2^2
		.
	\end{align*}
	It remains to show the estimate from below. To this end, we use that the function 
	\[
		f\FD{[1,\infty)}{[2,\infty)}6{x}{x+\frac 1x}
	\]
	is increasing, continuous, and convex. Its inverse function $f^{-1}$ thus exists and is concave, and we obtain
	\begin{align*}
		\|C\|^2 &\, =\|C^TC\|=
		f^{-1}\big(\Tr(C^TC)\big)
		\\
		= &\, f^{-1}\Big(
		\cos^2\sigma\!\cdot\! f\big(\max\big\{\big(\smfrac ab\big)^2,\big(\smfrac ba\big)^2\big\}\big)
		+\sin^2\sigma\!\cdot\! f\big(\max\big\{(ab)^2,\big(\smfrac 1{ab}\big)^2\big\}\big)
		\Big)
		\\
		\geq &\, 
		\cos^2\sigma\cdot\max\big\{\big(\smfrac ab\big)^2,\big(\smfrac ba\big)^2\big\}
		+\sin^2\sigma\cdot\max\big\{(ab)^2,\big(\smfrac 1{ab}\big)^2\big\}
		=\|v_+\|_2^2
	\end{align*}
\end{proof}

\noindent
The proof of the theorem is now easily completed. 

\begin{proof}[Proof of \Cref{Q6}]
	Given data as in the theorem, we apply \Cref{Q1} with the matrices 
	\[
		\Omega_j\DE\Omega(a_j,\psi_j),\quad j=1,\ldots,N
		,
	\]
	and use \Cref{Q5}. This yields 
	\begin{align*}
		\log\|W_H(z)\|\leq \mkern-100mu&\,\mkern100mu 
		|z|\sum_{j=1}^N \int_{y_{j-1}}^{y_j}\|\Omega(a_j,\psi_j)H(t)J\Omega(a_j,\psi_j)^{-1}\|\DD t
		\\
		&\, +\sum_{j=1}^{N-1}\log\|\Omega(a_j,\psi_j)\Omega(a_{j+1},\psi_{j+1})^{-1}\|
		\\
		&\, +\log\|\Omega(a_1,\psi_1)^{-1}\|+\log\|\Omega(a_N,\psi_N)\|
		\\
		\leq &\, |z|\sum_{j=1}^N\bigg(
		\int_{y_{j-1}}^{y_j}\Big(a_j^2\cos^2(\phi(t)-\psi_j)+\frac 1{a_j^2}\sin^2(\phi(t)-\psi_j)\Big)\cdot\Tr H(t)\DD t
		\bigg)
		\\
		&\, +\sum_{j=1}^{N-1}\log\bigg(
		\max\Big\{\frac{a_j}{a_{j+1}},\frac{a_{j+1}}{a_j}\Big\}\cdot\big|\cos\big(\psi_j-\psi_{j+1}\big)\big|
		+\frac{|\sin(\psi_j-\psi_{j+1})|}{a_ja_{j+1}}
		\bigg)
		\\
		&\, +\log\frac 1{a_1}+\log\frac 1{a_N}
		\\
		= &\, |z|(A_1+A_2)+A_3+A_4
		.
	\end{align*}
\end{proof}

\begin{Remark}
\label{Q11}
	The estimate stated in the theorem could be slightly improved on the cost of writing a much more cumbersome expression 
	$A_3'$ instead of $A_3$. Namely, by using the exact value for the norm in \Cref{Q5}\,\Enumref{3} instead of the upper
	estimate given there. Doing this would turn the inequality on the fourth line of the above estimate into an equality. 

	The upper and lower bounds for the norm in \Cref{Q5}\,\Enumref{3} differ only at most by the universal
	multiplicative constant $\sqrt 2$. Hence, the potential improvement is limited by 
	\[
		A_3\leq A_3'+(N-1)\cdot\frac 12\log 2
		.
	\]
\end{Remark}

\noindent
Let us now show that \cite[Theorem~1]{romanov:2017} can indeed be deduced from \Cref{Q6}. 
Recall the statement (for convenience we formulate Romanov's theorem in a notation already fitting \Cref{Q6}). 

\begin{theorem}[\cite{romanov:2017}]
\label{Q12}
	Let $H(t)=\xi_{\phi(t)}\xi_{\phi(t)}^T$ be a Hamiltonian on an interval $[0,L]$, and let $d\in(0,1)$. 
	Assume that we are given a constant $C>0$, and for each sufficiently large $R$ 
	\begin{Ilist}
	\item a partition $(y_0,\ldots,y_{N(R)})$ of $[0,L]$,
	\item rotation parameters $\psi_1(R),\ldots,\psi_{N(R)}(R)\in\bb R$,
	\item distortion parameters $a_1(R),\ldots,a_{N(R)}(R)\in(0,1]$,
	\end{Ilist}
	such that 
	\begin{Enumerate}
	\item ${\displaystyle
		\sum_{j=1}^{N(R)}\frac 1{a_j(R)^2}
		\int_{y_{j-1}(R)}^{y_j(R)}\big\|H(t)-\xi_{\psi_j(R)}\xi_{\psi_j(R)}^T\big\|\DD t
		\leq CR^{d-1}
		}$,
	\item ${\displaystyle
		\sum_{j=1}^{N(R)}a_j(R)^2\big(y_j(R)-y_{j-1}(R)\big)
		\leq CR^{d-1}
		}$,
	\item ${\displaystyle
		\sum_{j=1}^{N(R)-1}\log\Big(1+\frac{|\sin(\psi_j(R)-\psi_{j-1}(R))|}{a_j(R)a_{j+1}(R)}\Big)
		\leq CR^d
		}$,
	\item ${\displaystyle
		\log\frac 1{a_1(R)}+\log\frac 1{a_{N(R)}}+
		\sum_{j=1}^{N(R)-1}\Big|\log\frac{a_{j+1}(R)}{a_j(R)}\Big|
		\leq CR^{d-1}
		}$.
	\end{Enumerate}
	Then there exists a constant $K>0$ such that 
	\begin{equation}
	\label{Q71}
		\forall z\in\bb C\DP \log\|W_H(z)\|\leq K|z|^d
		.
	\end{equation}
\end{theorem}
\begin{proof}[Deduction from \Cref{Q6}]
	Let $R>0$ and assume that we have data $y_j(R),\psi_j(R),a_j(R)$ satisfying \Enumref{1}--\Enumref{4}.
	We are going to estimate the expressions $A_1,\ldots,A_4$ from \Cref{Q6}. 

	First, it is clear that 
	\begin{align*}
		A_1(R)= &\, \sum_{j=1}^{N(R)} a_j(R)^2\int_{y_{j-1}}^{y_j}\cos^2\big(\phi(t)-\psi_j(R)\big)\DD t
		\\
		\leq &\, \sum_{j=1}^{N(R)}a_j(R)^2\big(y_j(R)-y_{j-1}(R)\big)\leq CR^{d-1}
		.
	\end{align*}
	Next, observe that for all $\phi,\psi\in\bb R$ 
	\[
		\xi_\phi\xi_\phi^T-\xi_\psi\xi_\psi^T=\sin(\phi-\psi)\cdot
		\begin{pmatrix}
			-\sin(\phi+\psi) & \cos(\phi+\psi)
			\\
			\cos(\phi+\psi) & \sin(\phi+\psi)
		\end{pmatrix}
		.
	\]
	Since the matrix on the right side is unitary, it follows that 
	\[
		\|\xi_\phi\xi_\phi^T-\xi_\psi\xi_\psi^T\|=|\sin(\phi-\psi)|
		.
	\]
	From this we obtain 
	\begin{align}
		\label{Q72}
		A_2(R)= &\, 
		\sum_{j=1}^{N(R)}\frac 1{a_j(R)^2}\int_{y_{j-1}}^{y_j}\sin^2\big(\phi(t)-\psi_j(R)\big)\DD t
		\\
		\nonumber
		\leq &\, 
		\sum_{j=1}^{N(R)}\frac 1{a_j(R)^2}\int_{y_{j-1}}^{y_j}\big|\sin\big(\phi(t)-\psi_j(R)\big)\big|\DD t
		\leq CR^{d-1}
		.
	\end{align}
	Finally, we have 
	\begin{align*}
		A_3(R)+ &\, A_4(R)
		\\
		\leq &\, 
		\sum_{j=1}^{N(R)-1}\log\bigg[
		\max\Big\{\frac{a_j(R)}{a_{j+1}(R)},\frac{a_{j+1}(R)}{a_j(R)}\Big\}
		\bigg(1+\frac{|\sin(\psi_j-\psi_{j+1})|}{a_j(R)a_{j+1}(R)}\bigg)
		\bigg]
		\\
		&\, -\log a_1(R)-\log a_{N(R)}(R)
		\\
		= &\, 
		\sum_{j=1}^{N(R)-1}\Big|\log\frac{a_{j+1}(R)}{a_j(R)}\Big|+
		\sum_{j=1}^{N(R)-1}\Big(1+\frac{|\sin(\psi_j-\psi_{j+1})|}{a_j(R)a_{j+1}(R)}\Big)
		\\
		&\, +\log\frac 1{a_1(R)}+\log\frac 1{a_{N(R)}(R)}
		\leq 2CR^d
		.
	\end{align*}
	Now \cref{Q7} yields 
	\[
		\forall z\in\bb C\DP \log\|W_H(z)\|\leq|z|\cdot 2CR^{d-1}+2CR^d
		.
	\]
	We use this for $z\in\bb C$ with $|z|=R$ and obtain
	\[
		\log\|W_H(z)\|\leq 4CR^d\quad\text{for }|z|=R
		.
	\]
	By the assumption of the theorem the above argument can be made for all sufficiently large $R$. 
	Hence, \cref{Q71} follows.
\end{proof}

\begin{Remark}
\label{Q73}
	The improvement of \Cref{Q6} compared to \cite[Theorem~1]{romanov:2017} mainly happens in \cref{Q72}: clearly, 
	\[
		\sin^2(\phi(t)-\psi_j(R))\ll \big|\sin(\phi(t)-\psi_j(R)\big|
	\]
	when $\psi_j(R)$ is a good approximation of $\phi(t)$. 

	For this reason we also refer to \Cref{Q6} as the \emph{sine-square improvement} of Romanov's Theorem~1.
	We will see in \Cref{Q15} below that it is indeed a significant improvement.
\end{Remark}

\section{Hamiltonians with continuous rotation angle}
\label{Q76}

In this section we consider Hamiltonians of the form
\begin{equation}
\label{Q54}
	H(t)=\Tr H(t)\cdot\xi_{\phi(t)}\xi_{\phi(t)}^T,\quad t\in I
	,
\end{equation}
with a continuous rotation angle $\phi\DF I\to\bb R$, and prove an upper bound for $\log\|W_H(z)\|$. 
As a corollary we obtain a bound for the exponential order of the monodromy matrix of a H\"older continuous Hamiltonian which
improves \cite[Corollary~4\,(1.)]{romanov:2017}.
The proof of the upper estimate is an application of \Cref{Q6}, and nicely illustrates how concrete growth 
estimates can be deduced from the generic estimate. 

We use the following notation which involves the modulus of continuity of a function $\phi$. 
The case that $\phi$ is constant will be excluded, but this is no loss of generality: if $\phi$ in \cref{Q54} is constant, then 
$W$ is a linear polynomial and hence $\log\|W_H(z)\|=\BigO\big(\log|z|\big)$. 

\begin{Definition}
\label{Q53}
	For $\alpha\in\bb R\setminus\{0\}$ denote by $p_\alpha\DF(0,\infty)\to(0,\infty)$ the power function 
	$p_\alpha(t)\DE t^\alpha$. 
	\begin{Enumerate}
	\item Let $\omega\DF[0,\infty)\to[0,\infty)$ be a nondecreasing continuous function with $\omega(0)=0$ and 
		$\omega(\delta)>0$ for all $\delta>0$. Then we define an increasing bijection 
		$\Gamma_\omega\DF(0,\infty)\to(0,\infty)$ as
		\begin{equation}
		\label{Q27}
			\Gamma_\omega\DE p_{-1}\circ(p_1\cdot\omega)^{-1}\circ p_{-1}
			,
		\end{equation}
		where $(p_1\cdot\omega)^{-1}$ denotes the inverse function of $p_1\cdot\omega$
		(note here that $p_1\cdot\omega$ is an increasing bijection of $[0,\infty)$ onto itself).
	\item Let $H$ be a Hamiltonian of the form \cref{Q54} with continuous and non-constant rotation angle $\phi$. 
		Then we write $\omega_H$ for the modulus of uniform continuity of $\phi$, i.e.
		\[
			\omega_H(\delta)\DE\sup\big\{|\phi(t)-\phi(s)|\DS t,s\in I,|t-s|\leq\delta\big\}
			,\quad \delta\geq 0
			,
		\]
		and let $\Gamma_H\DE\Gamma_{\omega_H}$ be the function corresponding to $\omega_H$ by 
		the construction in item \Enumref{1}. 
	\end{Enumerate}
\end{Definition}

\noindent
The assignment $\omega\mapsto\Gamma_\omega$ defined by \cref{Q27} is injective. 
In fact, $\omega$ can be recovered from $\Gamma_\omega$ by the formula 
\[
	\omega=p_{-1}\cdot\big(p_{-1}\circ\Gamma_\omega^{-1}\circ p_{-1}\big)
	.
\]
Moreover, we have the following monotonicity property:
\[
	\omega_1\leq\omega_2\quad\Longrightarrow\quad\Gamma_{\omega_1}\leq\Gamma_{\omega_2}
	.
\]
Given a Hamiltonian $H$, the growth of the functions $\omega_H$ and $\Gamma_H$ is limited:
simply because $\omega_H$ is the modulus of uniform continuity of some continuous function on a compact interval, 
we have 
\[
	\omega_H(\delta)=\Smallo(1)\ \text{and}\ \delta=\BigO(\omega_H(\delta))\ \text{for}\ \delta\to 0
	.
\]
From this it follows that 
\begin{equation}
\label{Q55}
	\Gamma_H(r)=\Smallo(r)\ \text{and}\ \sqrt r=\BigO(\Gamma_H(r))\ \text{for}\ r\to\infty
	.
\end{equation}
Our bound for the monodromy matrix can now be formulated as follows. 

\begin{theorem}
\label{Q14}
	Let $H(t)=\Tr H(t)\cdot\xi_{\phi(t)}\xi_{\phi(t)}^T$ be a Hamiltonian on a compact interval $I=[\alpha,\beta]$ 
	whose rotation angle $\phi$ is continuous and not constant. Set $l\DE\beta-\alpha$ and $L\DE\int_I\Tr H(t)\DD t$. 
	Then 
	\begin{equation}
	\label{Q20}
		\log\|W_H(z)\|\leq 3l\cdot\Gamma_H\Big(\frac Ll|z|\Big)+\BigO\big(\log|z|\big)
		.
	\end{equation}
\end{theorem}
\begin{proof}
	We are going to apply \Cref{Q6}. 
	Let $\delta\in(0,l)$ and $a\in(0,1]$; a specific choice will be made later in dependence of $|z|$. 
	The data $y_j,\psi_j,a_j$ in \Cref{Q6} is now specified as follows:
	\begin{Ilist}
	\item Let $N$ be the unique positive integer with $N-1<\frac l\delta\leq N$, and define a partition 
		$(y_0,\ldots,y_N)$ of $I$ as
		\[
			y_j\DE
			\begin{cases}
				\alpha+j\cdot\delta &\text{if}\ j\in\{0,\ldots,N-1\}
				,
				\\
				\beta &\text{if}\ j=N
				.
			\end{cases}
		\]
	\item Rotation parameters are 
		\[
			\psi_j\DE\phi(y_j),\quad j=1,\ldots,N
			.
		\]
	\item Distortion parameters are $a_j\DE a$, $j=1,\ldots,N$.
	\end{Ilist}
	The choice of rotation parameters implies that 
	\[
		|\phi(t)-\psi_j|\leq\omega_H(\delta),\quad t\in[y_{j-1},y_j]
		.
	\]
	The constants $A_1,A_2,A_3$ from \Cref{Q6} can be estimated as follows:
	\begin{align*}
		A_1\leq &\, a^2L\ED B_1(a)
		,
		\\
		A_2\leq &\, \frac 1{a^2}\omega_H(\delta)^2L\ED B_2(a,\delta)
		,
		\\
		A_3\leq &\, 
		\sum_{j=1}^{N-1}\log\Big(1+\frac{|\sin(\psi_j-\psi_{j+1})|}{a^2}\Big)
		\\
		\leq &\, 
		(N-1)\frac 1{a^2}\omega_H(\delta)<\frac l{\delta}\frac 1{a^2}\omega_H(\delta)\ED B_3(a,\delta)
		.
	\end{align*}
	Given $z\in\bb C$ we specify the parameters $\delta$ and $a$ as
	\begin{equation}
	\label{Q23}
		\delta\DE\Big[\Gamma_H\Big(\frac{L}l|z|\Big)\Big]^{-1},\quad a\DE\omega_H(\delta)^{\frac 12}
		.
	\end{equation}
	These formulas are found by minimising the maximum of the expressions $B_1(a),B_2(a,\delta),B_3(a,\delta)$. 
	Observe that, by the properties of $\Gamma_H$ noted in \cref{Q55}, we have $\delta<l$ and $a\leq 1$ for 
	all sufficiently large $|z|$. We have 
	\begin{equation}
	\label{Q16}
		B_1(a)=B_2(a,\delta)=\omega_H(\delta)L,\quad B_3(a,\delta)=\frac l{\delta},\quad
		\delta\omega_H(\delta)=\frac l{L|z|}
		,
	\end{equation}
	and hence
	\[
		|z|B_2(a,\delta)=|z|B_1(a)=|z|\omega_H(\delta)L=\frac l{\delta}=B_3(a,\delta)
		.
	\]
	\Cref{Q6} implies that 
	\begin{align*}
		\log\|W_H(z)\|\leq &\, |z|\big(B_1(a)+B_2(a,\delta)\big)+B_3(a,\delta)+\log\frac 1{a^2}
		\\
		= &\, 3\cdot l\Gamma_H\Big(\frac{L}l|z|\Big)+\log\frac 1{\omega_H(\delta)}
		.
	\end{align*}
	By the last relation in \cref{Q16}
	\[
		\log\frac 1{\omega_H(\delta)}=\log\frac{L}l+\underbrace{\log\delta}_{<\log l}+\log|z|
		=\BigO\big(\log|z|\big)
		,
	\]
	and the bound \cref{Q20} follows. 
\end{proof}

\noindent
Applying \Cref{Q14} to H\"older continuous functions leads to the following corollary. 
To fix notation, recall that a function $\phi\DF I\to\bb R$ is called H\"older continuous with exponent $\alpha\in[0,1]$, if 
\begin{equation}
\label{Q74}
	\exists c>0\DQ\forall t,s\in I\DP |\phi(t)-\phi(s)|\leq c|t-s|^\alpha
	.
\end{equation}
The H\"older exponent $\alpha_1(\phi)$ of $\phi$ is 
\[
	\alpha_1(\phi)\DE\sup\big\{\alpha\in[0,1]\DS \phi\text{ is H\"older continuous with exponent }\alpha\big\}
	.
\]

\begin{corollary}
\label{Q18}
	Let $\alpha\in(0,1]$ and let $H(t)=\Tr H(t)\cdot\xi_{\phi(t)}\xi_{\phi(t)}^T$, $t\in I$, be a Hamiltonian on a 
	compact interval $I$ whose rotation angle is H\"older continuous with exponent $\alpha$.
	Then 
	\[	
		\log\Big(\max_{|z|=r}\|W_H(z)\|\Big)\lesssim r^{\frac 1{1+\alpha}}
		.
	\]
	Consequently, the exponential order of $W(z)$ does not exceed $\frac 1{1+\alpha_1(\phi)}$. 
\end{corollary}
\begin{proof}
	Let $c,\alpha$ be as in \cref{Q74} and set $\omega(\delta)\DE c\delta^\alpha$. 
	Then $\Gamma_\omega(r)=c^{\frac 1{1+\alpha}}r^{\frac 1{1+\alpha}}$. We have $\omega_H\leq\omega$, and thus also 
	$\Gamma_H\leq\Gamma_\omega$. \Cref{Q14} gives 
	\[
		\log\Big(\max_{|z|=r}\|W(z)\|\Big)\lesssim\Gamma_H\Big(\frac{L}lr\Big)+\BigO(\log r)
		\leq\Gamma_\omega\Big(\frac{L}lr\Big)+\BigO(\log r)\asymp r^{\frac 1{1+\alpha}}
		.
	\]
\end{proof}

\noindent
This corollary shows that the present generic estimate is an improvement of Romanov's Theorem even on the scale of
exponential order. 

\begin{Remark}
\label{Q15}
	In \cite[Corollary~4\,(1.)]{romanov:2017} it is shown that for a H\"older continuous 
	(trace normed) Hamiltonian with H\"older exponent $\alpha\in(0,1]$ the order of the entire function $W_H(z)$ does not 
	exceed $1-\frac\alpha 2$. The above corollary improves this: 	
	\[
		\forall \alpha\in(0,1)\DP \frac 1{1+\alpha}<1-\frac\alpha 2
		.
	\]
\end{Remark}

\noindent
\Cref{Q14} is limited to orders in $[\frac 12,1]$: due to \cref{Q55} the bound \cref{Q20} cannot go below $r^{\frac 12}$. 
To show that this really is a limitation, we should give an example of a Hamiltonian with continuous rotation angle and small 
order.

\begin{Example}
\label{Q26}
	We start from the example given in \cite[\S7.3]{romanov:2017}. Let $p\in(0,1)$, and let $\mu$ be a probability measure
	on $[0,1]$ which has no point masses, whose topological support has zero Lebesgue measure and is such that 
	its contiguous intervals, call them $I_j=(\alpha_j,\beta_j)$, satisfy 
	\[
		\sum_j(\beta_j-\alpha_j)^p<\infty
		.
	\]
	Set $\Delta\DE\bigcup_j\big(\mu([0,\alpha_j])+I_j\big)\subseteq[0,2]$, and note that $\Delta$ has measure $1$. 
	Let us show that $\Delta$ is dense in $[0,2]$. By our assumption that $\mu$ has no point masses, the function 
	$f(x)\DE\mu([0,x])+x$ is continuous. Given $t\in[0,2]$ we thus find $x\in[0,1]$ with $f(x)=t$. 
	The support of $\mu$ has empty interior, hence we can choose $x_n\in I_{j_n}$ such that $\lim_{n\to\infty}x_n=x$. 
	It follows that $t=f(x)=\lim_{n\to\infty}f(x_n)$, and since the distribution function of
	$\mu$ is constant on intervals $I_j$ 
	\[
		f(x_n)=\mu([0,x_n])+x_n=\mu([0,\alpha_{j_n}])+x_n\in\mu([0,\alpha_{j_n}])+I_{j_n}\subseteq\Delta
		.
	\]
	Let $H\DF[0,2]\to\bb R^{2\times 2}$ be the Hamiltonian defined as 
	\[
		H(t)\DE
		\begin{cases}
			\smmatrix 1000 &\text{if}\ t\in\Delta,
			\\[3mm]
			\smmatrix 0001 &\text{if}\ t\in[0,2]\setminus\Delta.
		\end{cases}
	\]
	We set $h_1=\mathds{1}_\Delta$ and $h_2=\mathds{1}_{\Delta^c}$, so that $H=\smmatrix{h_1}00{h_2}$.
	By \cite[\S7.3]{romanov:2017} (we write $\rho(\cdot)$ for the order of an entire function) 
	\[
		\rho(W_H)\leq\frac{2p}{p+1} 
		.
	\]
	Now we apply the general procedure \cite[Section~4]{kaltenbaeck.winkler.woracek:bimmel} to construct a non-diagonal 
	Hamiltonian. Set
	\[
		m_j(t)\DE\int_0^t h_j(s)\DD s,\ j=1,2
		.
	\]
	Then $m_j\DF[0,2]\to[0,1]$ are continuous, $m_2$ is a nondecreasing surjection and $m_1$ is an increasing bijection.
	This allows us to define a continuous Hamiltonian $\tilde H\DF[0,1]\to\bb R^{2\times 2}$ by 
	\[
		m\DE m_2\circ m_1^{-1},\qquad
		\tilde H\DE
		\begin{pmatrix}
			1 & -m
			\\
			-m & m^2
		\end{pmatrix}
		.
	\]
	By \cite[Lemma~4.1]{kaltenbaeck.winkler.woracek:bimmel} we have 
	\[
		(0,1)W_H(z)\binom 01=(0,1)W_{\tilde H}(z^2)\binom 01
		,
	\]
	and it follows that 
	\[
		\rho(W_{\tilde H})\leq\frac p{p+1}
		.
	\]
	Making an approriate choice of $p$, this becomes arbitrarily small. 

	Note that we can write $\tilde H$ in the form \cref{Q54} with the continuous rotation angle 
	\[
		\tilde\phi(x)\DE-\arctan m(x),\quad x\in[0,1]
		.
	\]
\end{Example}

\noindent
One interesting observation about the statement in \Cref{Q14} is that passing from $H$ to a reparameterisation does not change
the monodromy matrix, but may drastically change the modulus of continuity of the rotation angle and with it the bound on the
right side of \cref{Q20}. This fact can be used to improve the bound. 

Methodologically this is not a surprise; it reflects that in the proof of \Cref{Q14} we applied the generic estimate only with
equidistant partitions and the modulus of uniform continuity. Making a change of scale we can try to flatten out the rotation
angle on sections where it heavily oscillates, and by this make the quality of its continuity more even across 
the whole interval. 
The other way to achieve this effect would be to use arbitrary partitions. Maybe this would be more effective, but certainly it 
is computationally harder to handle. 

At this point let us just illustrate by an example that working with reparameterisations indeed can leed 
to a significant improvement. 

\begin{Example}
\label{Q19}
	For $\gamma,\beta>0$ let $\phi_{\gamma,\beta}\DF [0,1]\to\bb R$ be the chirp function 
	\[
		\phi_{\gamma,\beta}(t)\DE
		\begin{cases}
			t^\gamma\sin\big(\frac 1{t^\beta}\big) &\text{if}\ t\in(0,1]
			\\
			0 &\text{if}\ t=0
		\end{cases}
	\]
	and consider the Hamiltonian 
	\[
		H_{\gamma,\beta}(t)\DE \xi_{\phi_{\gamma,\beta}}\xi_{\phi_{\gamma,\beta}}^T,\quad t\in[0,1]
		.
	\]
	We require in the following that $\gamma\leq\beta$, so that $\phi_{\gamma,\beta}$ is not of bounded variation. 
	This is done to rule out an application of \cite[Corollary~4\,(2.)]{romanov:2017} which would imply at once that the 
	order of the monodromy matrix is at most $\frac 12$ (and we could not go below order $\frac 12$ anyway). 
	Our aim is to show that the order of the monodromy matrix $W_{H_{\gamma,\beta}}(z)$ is bounded by
	\begin{equation}
	\label{Q21}
		\rho\big(W_{H_{\gamma,\beta}}\big)\leq\frac{\beta}{\beta+\gamma}
		.
	\end{equation}
	The H\"older exponent of $\phi_{\gamma,\beta}$ is $\frac\gamma{\beta+1}$. Hence, \Cref{Q18} gives 
	\[
		\rho(W_{H_{\gamma,\beta}})\leq\frac 1{1+\frac\gamma{\beta+1}}=\frac{\beta+1}{\beta+\gamma+1}
		.
	\]
	For $\kappa>1$ set $\psi_\kappa(t)\DE t^\kappa$. Then $\psi_\kappa$ is an absolutely continuous increasing bijection of
	$[0,1]$ onto itself whose derivative is positive almost everywhere. It thus qualifies for being used as a
	reparameterisation. Denote 
	\[
		H_{\gamma,\beta}^{[\kappa]}(t)\DE \big(H_{\gamma,\beta}\circ\psi_\kappa\big)(t)\cdot\psi_\kappa'(t)
		,\quad t\in[0,1]
		.
	\]
	Apparently, $H_{\gamma,\beta}^{[\kappa]}=\kappa t^{\kappa-1}\cdot H_{\gamma\kappa,\beta\kappa}$, and hence 
	\[
		\rho\big(W_{H_{\gamma,\beta}}\big)=\rho\big(W_{H_{\gamma,\beta}^{[\kappa]}}\big)
		\leq\frac{\beta\kappa+1}{\beta\kappa+\gamma\kappa+1}
		=\frac{\beta+\frac1\kappa}{\beta+\gamma+\frac1\kappa}
		.
	\]
	Sending $\kappa$ to infinity, \cref{Q21} follows. 
\end{Example}

\section{Sharpness in \Cref{Q14}}

Remember \Cref{Q26} where the bound from \Cref{Q14} cannot possibly give the correct growth of the monodromy matrix. 
Our aim in this section is to construct examples where \cref{Q20} gives the correct growth, 
at least up to an error of logarithmic size.
In particular, in these examples, \cref{Q20} will give the correct order. We formulate this fact in a fairly general way. 

\begin{Theorem}
\label{Q49}
	Let $\ms g$ and $\ms m$ be regularly varying function with 
	\[
		\frac 12<\Ind\ms g<1\quad \text{and}\quad \int_1^\infty\frac 1{\ms m(t)}\DD t<\infty
		,
	\]
	and let $\ms n$ be regularly varying with $(\ms n\circ\ms m)(x)\sim(\ms m\circ\ms n)(x)\sim x$.

	Then there exists a Hamiltonian $H(t)=\xi_{\phi(t)}\xi_{\phi(t)}^T$ whose rotation angle $\phi(t)$ is continuous, such 
	that
	\begin{equation}
	\label{Q66}
		\big(\ms n\circ\ms g\big)(r)\lesssim\log\Big(\max_{|z|=r}\|W_H(z)\|\Big)\lesssim\ms g(r)
		.
	\end{equation}
\end{Theorem}

\noindent
Note that the gap left by \cref{Q66} is indeed rather small: we could choose for examples $\ms m(r)\DE r(\log r)(\log\log r)^2$. 
Then $\ms n(r)\sim\frac r{(\log r)(\log\log r)^2}$, and hence the lower bound \cref{Q66} satisfies
\[
	(\ms n\circ\ms g)(r)\asymp\frac{\ms g(r)}{(\log r)(\log\log r)^2}
	.
\]
In particular we see that in the H\"older continuous situation the bound for order given in \Cref{Q18} is sharp. 

For the proof of \Cref{Q49} we have to construct a function $\phi(t)$ whose modulus of continuity is prescribed and such that 
the growth of the corresponding monodromy matrix can be estimated from below. 

Finding just some function with given modulus of continuity is of course easy. Every continuous increasing and 
subadditiv function $\omega\DF[0,\infty)\to[0,\infty)$ with $\omega(0)=0$ is the modulus of continuity of itself. 
However, using such functions for the rotation angle $\phi(t)$ of a Hamiltonian will not lead to a required example:
the order of the monodromy matrix cannot exceed $\frac 12$ by \cite[Corollary~4(2.)]{romanov:2017}. 

It turns out that the following example of an oscillating function with prescribed modulus of continuity does the job. We want
to point out that placing constancy intervals is crucial, at least for our argument. 

\begin{Example}
\label{Q56}
	Assume we are given
	\begin{Enumerate}
	\item a sequence $(l_j)_{j=1}^\infty$ of positive numbers with $\sum_{j=1}^\infty l_j<\infty$, 
	\item a nonincreasing sequence $(m_j)_{j=1}^\infty$ of positive numbers with $m_j\leq l_j$ for all $j\in\bb N$,
	\item a continuous function $\pi\DF(0,\infty)\to(0,\infty)$, such that $\pi$ is nondecreasing on $(0,m_1)$, 
		the function $p_{-1}\cdot\pi$ (again $p_\alpha(x)\DE x^\alpha$) is nonincreasing on $(0,m_1)$, and 
		\[
			\lim_{x\to 0}\pi(x)=0,\quad \sup_{j\in\bb N}\frac{\pi(m_j)}{\pi(m_{j+1})}<\infty,\quad 
			\pi(m_1)<\frac\pi 2
			.
		\]
	\end{Enumerate}
	Set 
	\begin{align*}
		& r_n\DE\sum_{j=1}^{n-1}(l_j+m_j),\quad s_n\DE r_n+l_n,\quad L\DE\sum_{j=1}^\infty(l_j+m_j)
		,
		\\
		& \phi_n\DE\sum_{j=1}^{n-1}(-1)^{j+1}\pi(m_j)
		,
	\end{align*}
	and let $\phi\DF[0,L)\to[0,\frac\pi2)$ be the piecewise linear path connecting the points 
	\[
		(r_1,\,\phi_1),\,(s_1,\,\phi_1),\,(r_2,\,\phi_2),\,(s_2,\,\phi_2),\,(r_3,\,\phi_3),\,\ldots
	\]
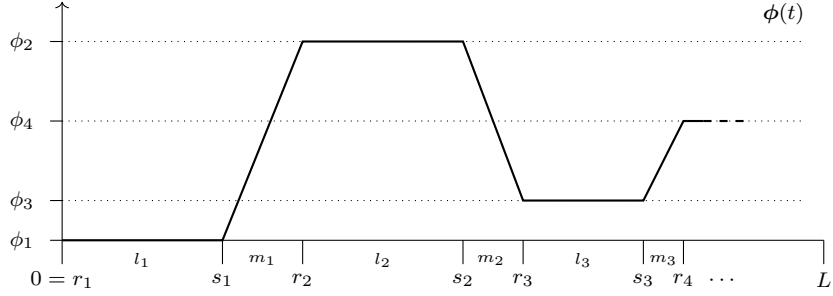
	\begin{center}
	\begin{tikzpicture}[x=1.5pt,y=1.5pt,scale=1,font=\fontsize{8}{8}]
		\draw[thin] (20,20)--(210,20);
		\draw[thin, ->] (20,20)--(20,80);

		\draw (10,20) node {$\phi_1$};
		\draw[thin] (17,20)--(20,20);
		\draw (10,70) node {$\phi_2$};
		\draw[thin] (17,70)--(20,70);
		\draw (10,30) node {$\phi_3$};
		\draw[thin] (17,30)--(20,30);
		\draw (10,50) node {$\phi_4$};
		\draw[thin] (17,50)--(20,50);

		\draw (20,10) node {$0=r_1$};
		\draw[thin] (20,14)--(20,20);
		\draw (60,10) node {$s_1$};
		\draw[thin] (60,14)--(60,20);
		\draw (80,10) node {$r_2$};
		\draw[thin] (80,14)--(80,20);
		\draw (120,10) node {$s_2$};
		\draw[thin] (120,14)--(120,20);
		\draw (135,10) node {$r_3$};
		\draw[thin] (135,14)--(135,20);
		\draw (165,10) node {$s_3$};
		\draw[thin] (165,14)--(165,20);
		\draw (175,10) node {$r_4$};
		\draw[thin] (175,14)--(175,20);
		\draw (185,10) node {$\cdots$};
		\draw (210,10) node {$L$};
		\draw[thin] (210,14)--(210,20);
		
		\draw (40,15) node {${\scriptstyle l_1}$};
		\draw (70,15) node {${\scriptstyle m_1}$};
		\draw (100,15) node {${\scriptstyle l_2}$};
		\draw (127,15) node {${\scriptstyle m_2}$};
		\draw (150,15) node {${\scriptstyle l_3}$};
		\draw (170,15) node {${\scriptstyle m_3}$};

		\draw (200,77) node {$\bm\phi(t)$};

		\draw[thick] (20,20)--(60,20)--(80,70)--(120,70)--(135,30)--(165,30)--(175,50)--(180,50);
		\draw[thick, dashed] (180,50)--(190,50);
		\draw[thin, dotted] (20,30)--(205,30);
		\draw[thin, dotted] (20,50)--(205,50);
		\draw[thin, dotted] (20,70)--(205,70);
	\end{tikzpicture}
	\end{center}
	We assert that the modulus of uniform continuity $\omega_\phi$ of the function $\phi$ satisfies 
	\begin{equation}
	\label{Q67}
		\omega_\phi(\delta)\asymp\pi(\delta)\quad\text{for }\delta\leq m_1
		.
	\end{equation}
	In order to prove this, we first show that 
	\begin{equation}
	\label{Q68}
		\omega_\phi(m_n)=\pi(m_n)
		.
	\end{equation}
	The inequality ``$\geq$'' follows since we have
	\[
		r_{n+1}-s_n=m_n\text{ and }|\phi(r_{n+1})-\phi(s_n)|=|\phi_{n+1}-\phi_n|=\pi(m_n)
		.
	\]
	Consider two points $t,s$ with $t<r_n$ and $t\leq s\leq t+m_n$. Then $s\leq s_n$. 
	The function $\phi|_{[0,s_n]}$ is a polygonal path with maximal slope $\frac{\pi(m_{n-1})}{m_{n-1}}$, and we obtain 
	\[
		|\phi(t)-\phi(s)|\leq\frac{\pi(m_{n-1})}{m_{n-1}}\cdot m_n\leq\pi(m_n)
		.
	\]
	For each two points $t,s$ with $r_n\leq t\leq s\leq t+m_n$ we have 
	\[
		|\phi(t)-\phi(s)|\leq|\phi_{n+1}-\phi_n|=\pi(m_n)
		,
	\]
	and ``$\leq$'' in \cref{Q68} follows. 

	For the proof of \cref{Q67}, let $\delta\leq m_1$ be given. Let $n\in\bb N$ be such that 
	$m_{n+1}<\delta\leq m_n$, then 
	\[
		\frac{\pi(m_{n+1})}{\pi(m_n)}=\frac{\pi(m_{n+1})}{\omega_\phi(m_n)}\leq\frac{\pi(\delta)}{\omega_\phi(\delta)}
		\leq\frac{\pi(m_n)}{\omega_\phi(m_{n+1})}=\frac{\pi(m_n)}{\pi(m_{n+1})}
		.
	\]
\end{Example}

\begin{proof}[Proof of \Cref{Q49}]
	Based on \cite[Theorems~1.8.2,1.8.5]{bingham.goldie.teugels:1989} we may assume w.l.o.g.\ that $\ms g$ and $\ms m$ are 
	increasing bijections of $(0,\infty)$ onto itself. Moreover, we may say that $\ms m(1)$ is as large as it pleases us
	(and a concrete request will be put later). 

	We use \Cref{Q56} with the data 
	\[
		m_j=l_j\DE\frac 1{\ms m(j)},\quad \pi\DE p_{-1}\cdot\big(p_{-1}\circ\ms g^{-1}\circ p_{-1}\big)
		.
	\]
	We have to check that the conditions required in \Cref{Q56}\,(i)--(iii) are fullfilled. First, 
	\[
		\sum_{j=1}^\infty l_j\leq l_1+\int_1^\infty\frac 1{\ms m(t)}\DD t<\infty
		,
	\]
	and $m_j$ ($=l_j$) is decreasing. Second, $\pi$ is continuous and regularly varying (at $0$) with 
	\[
		\Ind\pi=\frac 1{\Ind\ms g}-1\in(0,1)
		.
	\]
	Hence, sufficiently close to $0$, $\pi$ is increasing and $p_{-1}\circ\pi$ is decreasing. 
	Now we assume (w.l.o.g.) that $\ms m(1)$ is so large that $\frac 1{\ms m(1)}$ is already sufficiently close to $0$ in
	the above sense and $<\frac\pi 2$. 
	Also the function $\pi\circ p_{-1}\circ\ms m$ is regularly varying, and hence 
	\[
		\lim_{j\to\infty}\frac{\pi(m_j)}{\pi(m_{j+1})}
		=\lim_{j\to\infty}\frac{(\pi\circ p_{-1}\circ\ms m)(j)}{(\pi\circ p_{-1}\circ\ms m)(j+1)}=1
		.
	\]
	In particular, the quotient is bounded. 

	Let $\phi$ be the function constructed in \Cref{Q56}. Then $\omega_\phi\asymp\pi$, and since $\pi$ is regularly varying
	it follows that $\Gamma_{\omega_\phi}\asymp\Gamma_\pi$. The latter function computes as 
	\[
		\Gamma_\pi=p_{-1}\circ(p_1\cdot\pi)^{-1}\circ p_{-1}=\ms g
		.
	\]
	Let $H$ be the Hamiltonian $H(t)\DE\xi_{\phi(t)}\xi_{\phi(t)}^T$. The upper bound in \cref{Q66} is just \cref{Q20}. 
	In order to show the lower bound, we aim at an application of \Cref{Q17,Q48}. We use the set 
	\[
		\Delta\DE\bigcup_{j=1}^\infty(s_j,r_j)
	\]
	in \Cref{Q17}. The Hamiltonian $\tilde H$ constructed there is in our situation the Hamburger Hamiltonian with lengths 
	$(l_j)_{j=1}^\infty$ and angles $(\phi_j)_{j=1}^\infty$. Now \Cref{Q48} comes into play: we have 
	\[
		l_{j+1}l_j\sin^2(\phi_{j+1}-\phi_j)=\frac 1{\ms m(j+1)\ms m(j)}\sin^2(\pi(m_j))
		\sim\frac 1{\ms m(j)^2}(\pi\circ p_{-1}\circ\ms m)(j)^2
		,
	\]
	and hence we can use 
	\[
		\ms f\DE p_{-2}\circ(p_1\cdot\pi)\circ p_{-1}\circ\ms m
	\]
	in \cref{Q69}. The right side of \cref{Q84} then is 
	\[
		\ms f^{-1}\circ p_2=\ms m^{-1}\circ p_{-1}\circ(p_1\cdot\pi)^{-1}\circ p_{-\frac 12}\circ p_2
		=\ms m^{-1}\circ\ms g
		.
	\]
	Thus, by \Cref{Q48}, 
	\[
		\log\max_{|z|=r}\|W_{\tilde H}(z)\|\gtrsim(\ms m^{-1}\circ\ms g)(r)
		.
	\]
	\Cref{Q17} implies that also 
	\[
		\log\max_{|z|=r}\|W_H(z)\|\gtrsim(\ms m^{-1}\circ\ms g)(r)
		.
	\]
	Note here that $\Ind\ms m\geq 1$, and hence $\Ind(\ms m^{-1}\circ\ms g)\in(0,1)$. 
\end{proof}


{\footnotesize
\begin{flushleft}
	R.~Pruckner \\
	Institute for Analysis and Scientific Computing \\
	Vienna University of Technology \\
	Wiedner Hauptstra{\ss}e 8--10/101 \\
	1040 Wien \\
	AUSTRIA \\
	email: raphael.pruckner@tuwien.ac.at \\[5mm]
\end{flushleft}
\begin{flushleft}
	H.\,Woracek\\
	Institute for Analysis and Scientific Computing\\
	Vienna University of Technology\\
	Wiedner Hauptstra{\ss}e\ 8--10/101\\
	1040 Wien\\
	AUSTRIA\\
	email: \texttt{harald.woracek@tuwien.ac.at}\\[5mm]
\end{flushleft}
}


\begin{thebibliography}{99}
%
\bibitem{behrndt.hassi.snoo:2020}
J. Behrndt, S. Hassi, and H. de Snoo, \emph{Boundary value problems,
Weyl functions, and differential operators}, vol. 108, Monographs in
Mathematics, Birkh\"auser/Springer, Cham, 2020.%
\bibitem{berg.szwarc:2014}
C. Berg and R. Szwarc, “On the order of indeterminate moment
problems”, \emph{Adv. Math.} 250 (2014), pp. 105–143.
%
\bibitem{bingham.goldie.teugels:1989}
N.H. Bingham, C.M. Goldie, and J.L. Teugels, \emph{Regular
variation}, vol. 27, Encyclopedia of Mathematics and its Applications,
Cambridge University Press, Cambridge, 1989.
%
\bibitem{debranges:1968}
L. de Branges, \emph{Hilbert spaces of entire functions}, Englewood
Cliffs, N.J.: Prentice-Hall Inc., 1968.
%
\bibitem{dunford.schwartz:1963}
N. Dunford and J. T. Schwartz, \emph{Linear operators. Part II:
Spectral theory. Self adjoint operators in Hilbert space}, With the
assistance of William G. Bade and Robert G. Bartle, Interscience
Publishers John Wiley \& Sons New York-London, 1963.
%
\bibitem{hassi.snoo.winkler:2000}
S. Hassi, H. de Snoo, and H. Winkler, “Boundary-value problems for
two-dimensional canonical systems”, \emph{Integral Equations Operator
Theory} 36.4 (2000), pp. 445–479.
%
\bibitem{kac:1984}
I.S. Kac, “Linear relations, generated by a canonical differential equation
on an interval with a regular endpoint, and expansibility in
eigenfunctions”, Russian, \emph{VINITI Deponirovannye Nauchnye
Raboty} 195.1 (1985), Deposited in Ukr NIINTI, No. 1453, 1984, 50 pp.,
b.o. 720.
%
\bibitem{kac:1986a}
I.S. Kac, “Expansibility in eigenfunctions of a canonical differential
equation on an interval with singular endpoints and associated linear
relations”, Russian, \emph{VINITI Deponirovannye Nauchnye Raboty}
282.12 (1986), Deposited in Ukr NIINTI, No. 2111, 1986, 64 pp., b.o.
1536.
%
\bibitem{kac:1999}
I.S. Kac, “Inclusion of the Hamburger power moment problem in the
spectral theory of canonical systems”, Russian, \emph{Zap. Nauchn.
Sem. S.-Peterburg. Otdel. Mat. Inst. Steklov. (POMI)} 262.Issled. po
Linein. Oper. i Teor. Funkts. 27 (1999), English translation: J. Math.
Sci. (New York) 110 (2002), no. 5, 2991–3004, pp. 147–171, 234.
%
\bibitem{kac:2002}
I.S. Kac, “Linear relations generated by a canonical differential equation
of dimension 2, and eigenfunction expansions”, \emph{Algebra i Analiz}
14.3 (2002), pp. 86–120.
%
\bibitem{kaltenbaeck.winkler.woracek:bimmel}
M. Kaltenb\"ack, H. Winkler, and H. Woracek, “Strings, dual strings, and
related canonical systems”, \emph{Math. Nachr.} 280.13-14 (2007),
pp. 1518–1536.
%
\bibitem{levin:1980}
B.Ja. Levin, \emph{Distribution of zeros of entire functions}, Revised,
vol. 5, Translations of Mathematical Monographs, Translated from the
Russian by R. P. Boas, J. M. Danskin, F. M. Goodspeed, J. Korevaar, A.
L. Shields and H. P. Thielman, Providence, R.I.: American Mathematical
Society, 1980.
%
\bibitem{livshits:1939}
M.S. Liv\v sic, “On some questions concerning the determinate case of
Hamburger’s moment problem”, Russian. English summary,
\emph{Rec. Math. N. S. [Mat. Sbornik]} 6(48) (1939), pp. 293–306.
%
\bibitem{orcutt:1969}
B.C. Orcutt, \emph{Canonical differential equations}, Thesis
(Ph.D.)–University of Virginia, ProQuest LLC, Ann Arbor, MI, 1969.
%
\bibitem{pruckner.romanov.woracek:jaco-ASC}
R. Pruckner, R. Romanov, and H. Woracek, \emph{Bounds on order of
indeterminate moment sequences (extended preprint)}, 26 pp., ASC
Report 38,
http://www.asc.tuwien.ac.at/preprint/2015/asc38x2015.pdf, Vienna
University of Technology, 2015.
%
\bibitem{pruckner.romanov.woracek:jaco}
R. Pruckner, R. Romanov, and H. Woracek, “Bounds on order of
indeterminate moment sequences”, \emph{Constr. Approx.} 46 (2017),
pp. 199–225.
%
\bibitem{pruckner.woracek:srt}
R. Pruckner and H. Woracek, “Estimates for the order of Nevanlinna
matrices and a Berezanskii-type theorem”, \emph{Proc. Roy. Soc.
Edinburgh Sect. A} 149.6 (2019), pp. 1637–1661.
%
\bibitem{remling:2018}C. Remling, \emph{Spectral Theory of Canonical Systems}, De Gruyter
Studies in Mathematics Series, Walter de Gruyter GmbH, 2018.
%
\bibitem{romanov:1408.6022v1}
R. Romanov, \emph{Canonical systems and de Branges spaces},
version 1, Aug. 26, 2014, arXiv: 1408.6022v1[math.SP].
%
\bibitem{romanov:2017}
R. Romanov, “Order problem for canonical systems and a conjecture of
Valent”, \emph{Trans. Amer. Math. Soc.} 369.2 (2017),
pp. 1061–1078.
\end{thebibliography}
\end{document}